\documentstyle[twoside,11pt]{article}

\newtheorem{theorem}{Theorem}[section]
\newtheorem{lemma}[theorem]{Lemma}
\newtheorem{corollary}[theorem]{Corollary}
\newtheorem{statement}[theorem]{Statement}

\long\gdef\boxit#1{\begingroup\vbox{\hrule\hbox{\vrule\kern3pt
      \vbox{\kern3pt#1\kern3pt}\kern3pt\vrule}\hrule}\endgroup}

\def\qed{ \ \vrule width.2cm height.2cm depth0cm\smallskip}

\def\Keywords{\bigskip\par {\sl Keywords\/}:\enspace}

\def\max{\,{\rm max}}
\def\min{\,{\rm min}}
\def\exten{{\rm Ext}}

\def\Rset{\mbox{\bf R}}
\def\Zset{\mbox{\bf Z}}

\def\bar{\overline}
\def\eps{\epsilon}
\def\taustar{\tau^\ast}
\def\kott{K^-_{3,3}}

\def\isom{\simeq}

\def\Xscr{{\cal X}}
\def\Fscr{{\cal F}}

\def\Vscr{{\cal V}}

\def\Escr{{\cal E}}
\def\Ascr{{\cal A}}
\def\Bscr{{\cal B}}

\def\Eie{E_i^=}
\def\Eine{E_i^{\ne}}

\def\Eonene{E_1^{\ne}}

\def\Etwone{E_2^{\ne}}
\newenvironment{proof}{\noindent {\bf Proof.\/}}{$\qed$\vskip 0.1in}
\newcommand{\Xcomment}[1]{}

\newenvironment{myitem}{\refstepcounter{equation}\begin{enumerate}%
\item[(\thesection.\arabic{equation})]$\quad$}{\end{enumerate}}

\makeatletter
\@addtoreset{equation}{section}

\makeatother

\begin{document}
\begin{titlepage}
\title{\bf How to uncross some modular metrics
\thanks{Supported by the Fields Institute for Research in Mathematical
Sciences and by grant 97-01-00115 from the Russian Foundation
of Basic Research.}
}

\def\thepage {} 

\author{{\Large Alexander V. Karzanov}
\thanks{On leave from the
Institute for System Analysis of the Russian Acad. Sci.,
Prospect 60 Let Oktyabrya, 117312 Moscow, Russia,
email: {\sl sasha@cs.isa.ac.ru}
}
}

\date{December 1999}

\maketitle

\begin{center}
Dedicated to the memory of Walter Deuber
\end{center}

\bigskip

  \begin{abstract}
Let $\mu$ be a metric on a set $T$, and let $c$ be a nonnegative
function on the unordered pairs of elements of a superset
$V\supseteq T$. We consider the problem of minimizing the inner
product $c\cdot m$
over all semimetrics $m$ on $V$ such that $m$ coincides with $\mu$
within $T$ and each element of $V$ is at zero distance from $T$
(a variant of the {\em multifacility location problem}). In
particular, this generalizes the well-known multiterminal (or
multiway) cut problem.

Two cases of metrics $\mu$ have been known for which the problem can
be solved in polynomial time: (a) $\mu$ is a modular metric whose
underlying graph $H(\mu)$ is hereditary modular and orientable
(in a certain sense); and (b) $\mu$ is a median metric. In the
latter case an optimal solution can be found by use of
a cut uncrossing method.

\Xcomment{
We give a common generalization for both cases by proving that the
problem is in P for any modular metric $\mu$ whose all orbit graphs are
hereditary modular and orientable. To this aim, we show the existence of
a retraction of the Cartesian product of the orbit graphs to $H(\mu)$,
which enables us to elaborate an analog of the cut uncrossing method
for such metrics $\mu$.
}
In this paper we generalize the idea of cut uncrossing to show the
polynomial solvability for a wider class of metrics $\mu$, which
includes the median metrics as a special case. The metric
uncrossing method that we develop relies on the existence of
retractions of certain modular graphs. On the negative side,
we prove that for $\mu$ fixed, the problem is NP-hard if $\mu$ is
non-modular or $H(\mu)$ is non-orientable.
  \end{abstract}

\Keywords location problem, multiterminal (multiway) cut,
modular graph


  \end{titlepage}

\baselineskip 15pt

\section{Introduction} \label{sec:introduc}

We deal with a variant of the {\em multifacility location problem}.
In its setting, there are a finite metric space $(T,\mu)$, a finite
set $X$,
and a nonnegative function $c$ on the pairs of elements of $T\cup X$.
(The elements of $T$ are thought of as the points where the existing
facilities are located, the elements of $X$ as new facilities,
and $c(x,y)$ as a measure of communication between $x$ to $y$.)
The objective is to place each new facility at a point of $T$ minimizing
the sum of values $c(x,y)\mu(x',y')$, where $x,y$ range over the pairs
of facilities and $x',y'$ are the points of $T$ where $x,y$ are
placed. For a survey on location problems, see, e.g.,~\cite{TFL-83}.

This problem can be reformulated in
terms of metric extensions. We start with some terminology and
notation. A {\it semimetric} on a set $S$ is a function
$d:S\times S\to\Rset_+$ that establishes {\it distances} on
the pairs of elements ({\em points})
of $S$ satisfying $d(x,x)=0$, $d(x,y)=d(y,x)$ and $d(x,y)+d(y,z)
\ge d(x,z)$, for all $x,y,z\in S$. We use notation $xy$ for an
unordered pair $\{x,y\}$ and usually write $d(xy)$ instead of
$d(x,y)$. The set of pairs $xy$ with $x\ne y$ is denoted by $E_S$.
When $d(xy)>0$ for all $xy\in E_S$, $d$ is a {\em metric}. We do not
distinguish between the (semi)metric $d$ and the (semi)metric space
$(S,d)$ and usually deal with only finite (semi)metric spaces.
A special case is the {\em path metric}\/ $d^G$ of a connected
graph $G$, where $d^G(xy)$ is the minimum number of edges
of a path in $G$ connecting nodes $x$ and $y$.

A semimetric $m$ on a set $V\supseteq S$ is said to be
an {\em extension} of $d$
if the restriction ({\em submetric}) of $m$ to $S$ is just $d$.
Such an $m$ is called a 0-{\em extension}\/ if the distance
$m(x,S)$ from each point $x\in V$ to $S$ is zero, i.e.,
$m(xs)=0$ for some $s\in S$. Clearly each 0-extension $m$ is
uniquely determined by the 0-{\em distance sets}\/
$X_s=\{x\in V: m(xs)=0\}$, $s\in S$, and these sets give a
partition of $V$ when $d$ is a metric.

The above problem is equivalent to the {\em minimum
0-extension problem} : Given a metric $\mu$ on a set $T$, a superset
$V\supseteq T$, and a function $c:E_V\to \Zset_+$,
   \begin{myitem}
{\sl Find a 0-extension $m$ of $\mu$ to $V$ with
$c\cdot m:=\sum(c(e)m(e): e\in E_V)$ minimum}.
  \label{eq:0ext}
  \end{myitem}

In this paper we extend earlier results on
the complexity of~(\ref{eq:0ext}) for fixed metrics $\mu$.

Two classes of metrics $\mu$ have been found for
which~(\ref{eq:0ext}) is solvable in polynomial time.
One class consists of the metrics for which~(\ref{eq:0ext})
becomes as easy as its linear programming
relaxation. More precisely, let $\tau=\tau(V,c,\mu)$ denote the
minimum $c\cdot m$ in~(\ref{eq:0ext}), and let
$\taustar=\taustar(V,c,\mu)$ denote the
minimum $c\cdot m$ in the problem:
  \begin{myitem}
{\sl Find an extension $m$ of $\mu$ to $V$ with
$c\cdot m$ minimum}.
 \label{eq:ext}
  \end{myitem}
Then $\tau\ge\taustar$. A metric $\mu$ is called {\em minimizable}\/
if $\tau(V,c,\mu)=\taustar(V,c,\mu)$
holds for any $V$ and $c$. Since (\ref{eq:ext}) is a linear program
whose constraint matrix size is polynomial in $|V|$, (\ref{eq:ext}) is
solvable in strongly polynomial time. This easily implies that for every
minimizable metric $\mu$, on optimal 0-extension in~(\ref{eq:0ext})
can be found in
strongly polynomial time as well. The following theorem characterizes
the set of minimizable path metrics.

\begin{theorem} \label{tm:mgraph} {\rm\cite{kar-98}}
For a graph $H$, the metric $d^H$ is minimizable if and only if $H$ is
hereditary modular and orientable.
 \end{theorem}

\begin{figure}[tb]
  \unitlength=1mm
 \begin{center}
  \begin{picture}(125,35)
\put(17,5){\circle*{1.5}}  
\put(17,27){\circle*{1.5}}
\put(39,5){\circle*{1.5}}
\put(39,27){\circle*{1.5}}
\put(12,3){$v_0$}
\put(12,27){$v_1$}
\put(41,27){$v_2$}
\put(41,3){$v_3$}
\put(17,5){\line(0,1){22}}
\put(17,5){\line(1,0){22}}
\put(17,27){\line(1,0){22}}
\put(39,5){\line(0,1){22}}
\put(35.5,4){$>$}
\put(35.5,26.1){$>$}
\put(15.75,24){$\wedge$}
\put(37.75,24){$\wedge$}
\put(90,5){\circle*{1.5}}  
\put(110,5){\circle*{1.5}}
\put(80,17.5){\circle*{1.5}}
\put(120,17.5){\circle*{1.5}}
\put(90,30){\circle*{1.5}}
\put(110,30){\circle*{1.5}}
\put(90,5){\line(1,0){20}}
\put(90,30){\line(1,0){20}}
\put(90,5){\line(-4,5){10}}
\put(110,5){\line(4,5){10}}
\put(90,30){\line(-4,-5){10}}
\put(110,30){\line(4,-5){10}}
\put(90,5){\line(4,5){20}}
\put(110,5){\line(-4,5){20}}
  \end{picture}
 \end{center}
\caption{\hspace{1cm} (a) An orientation of a 4-circuit \hspace{2cm}
  (b) $\kott$ \hspace{3cm} } \label{fig:orient}
  \end{figure}
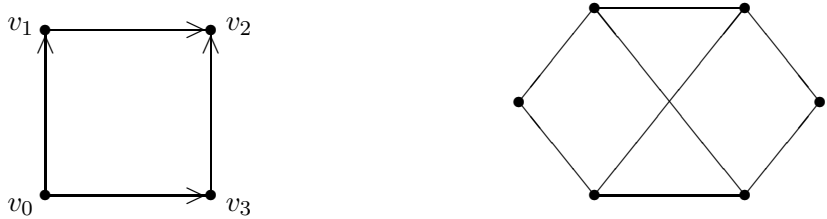

Recall that a metric $\mu$ on $T$ is {\em modular}\/ if every
three points $s_0,s_1,s_2\in T$ have a {\em median}, a node $z\in T$
satisfying $\mu(s_iz)+\mu(zs_j)=\mu(s_is_j)$ for all $0\le i<j\le
2$. A graph $H$ is called {\em modular}\/ if $d^H$ is
modular, and {\em hereditary modular}\/ if every isometric subgraph of
$H$ is modular, where a subgraph (or circuit) $H'=(T',U')$ of $H$ is
{\em isometric}\/ if $d^{H'}(st)=d^H(st)$ for all $s,t\in T'$.
Every modular graph is bipartite. A graph is called {\em
orientable}\/ if its edges
can be oriented so that for any 4-circuit $C=(v_0,e_1,v_1,...,e_4,
v_4=v_0)$ and $i=1,2$, the edge $e_i$ is oriented from
$v_{i-1}$ to $v_i$ if and
only if the opposite edge $e_{i+2}$ is oriented from $v_{i+2}$ to
$v_{i+1}$; see Fig.~\ref{fig:orient}(a). For example, every bipartite
graph with at most five nodes is hereditary modular and
orientable. The simplest hereditary modular but not orientable graph
is the graph $\kott$ obtained from $K_{3,3}$ by deleting one edge; see
Fig.~\ref{fig:orient}(b). Using terminology in~\cite{kar-98},
we refer to an orientable hereditary modular graph as a {\em
frame}.

Theorem~\ref{tm:mgraph} is extended to general metrics using
the notion of {\em underlying graph} of $\mu$. This is the least graph
$H(\mu)=(T,U(\mu))$ which enables us to restore $\mu$ if we know
the distances of its edges. Formally, nodes $x,y\in T$ are adjacent in
$H(\mu)$ if and only if no other node $z\in T$ lies {\em between}\/
$x$ and $y$, i.e., satisfies $\mu(xz)+\mu(zy)=\mu(xy)$.
This graph is modular if $\mu$ is modular~\cite{ban-85}.

 \begin{theorem} \label{tm:mmetr} {\rm \cite{BCK-99}}
A metric $\mu$ is minimizable if and only if $\mu$ is modular and
$H(\mu)$ is a frame.
  \end{theorem}

Another tractable case involves {\em median metrics}, the metrics
$\mu$ with precisely one median for each triplet of points.
Chepoi~\cite{che-97} showed that~(\ref{eq:0ext}) with any
median metric $\mu$ is solvable in strongly polynomial time.
A simple alternative method, based on cut uncrossing techniques, is
suggested in~\cite{kar-98}.
Note that a minimizable metric need not be a median one, and vice
versa. For example, $d^{K_{2,3}}$ is minimizable but not median, while
the path metric of the (skeleton of the) cube is median but not
minimizable (the cube is not hereditary modular as it contains an
isometric 6-circuit).

In this paper we show the polynomial solvability for
a class of modular metrics which includes the median ones as
a special case. It uses the notion of orbit graphs that we now
introduce. Given a modular graph $H=(T,U)$,
two edges are called {\em mates}\/ if they are opposite in some
4-circuit; when dealing with graphs with possible parallel edges,
we refer to such edges as mates as well. Edges $e,e'$ of $H$ are
called {\em projective}\/ if there is a sequence
$e=e_0,e_1,\ldots,e_k=e'$ of edges such that any two consecutive
$e_i,e_{i+1}$ are mates; such a sequence is called {\em
projective}\/ too. A maximal set $Q$ of mutually projective edges
is called an {\em orbit}. Define the {\em orbit graph}\/ $H_Q$
to be $H//(U-Q)$, where for a graph $H'$ and a subset $Z$ of its
edges, $H'//Z$ denotes the graph obtained by contracting $Z$ (i.e.,
forming $H'/Z$) and then identifying parallel edges appeared.

The main result in this paper is the following.

\Xcomment{
 \begin{theorem} \label{tm:main}
Let $\mu$ be a modular metric such that for each orbit $Q$ of the
underlying graph $H$ of $\mu$, the orbit graph $H_Q$ is a frame.
Then~(\ref{eq:0ext}) can be solved in strongly polynomial time.
  \end{theorem}
}
 \begin{theorem} \label{tm:main}
Let $\mu$ be a modular metric with underlying graph $H=(T,U)$, and
let for each orbit $Q$ of $H$,

{\rm (i)} the orbit graph $H_Q$ is a frame, and

{\rm (ii)} $H_Q$ is isomorphic to some subgraph of the graph
$(T,Q)$.

\noindent Then~(\ref{eq:0ext}) can be solved in strongly polynomial
time.
  \end{theorem}
%
%
\begin{figure}[tb]
 \unitlength=1.0mm
\begin{center}
  \begin{picture}(150,30)
\put(0,0){\circle*{1.5}}   
\put(20,0){\circle*{1.5}}
\put(40,0){\circle*{1.5}}
\put(0,20){\circle*{1.5}}
\put(20,20){\circle*{1.5}}
\put(40,20){\circle*{1.5}}
\put(10,10){\circle*{1.5}}
\put(48,10){\circle*{1.5}}
\put(68,10){\circle*{1.5}}
\put(88,10){\circle*{1.5}}
\put(48,30){\circle*{1.5}}
\put(68,30){\circle*{1.5}}
\put(88,30){\circle*{1.5}}
\put(58,20){\circle*{1.5}}
\put(115,10){\circle*{1.5}}   
\put(135,10){\circle*{1.5}}
\put(115,30){\circle*{1.5}}
\put(135,30){\circle*{1.5}}
\put(125,20){\circle*{1.5}}
\put(115,0){\circle*{1.5}}
\put(135,0){\circle*{1.5}}
\put(115,5){\circle*{1.5}}
\put(135,5){\circle*{1.5}}
\put(0,0){\line(1,0){20}}   
\put(0,20){\line(1,0){20}}
\put(48,10){\line(1,0){20}}
\put(48,30){\line(1,0){20}}
\put(115,10){\line(1,0){20}}
\put(115,30){\line(1,0){20}}
\put(0,0){\line(0,1){20}}
\put(20,0){\line(0,1){20}}
\put(40,0){\line(0,1){20}}
\put(48,10){\line(0,1){20}}
\put(68,10){\line(0,1){20}}
\put(88,10){\line(0,1){20}}
\put(115,10){\line(0,1){20}}
\put(135,10){\line(0,1){20}}
\put(0,20){\line(1,-1){20}}
\put(48,30){\line(1,-1){20}}
\put(115,30){\line(1,-1){20}}
\linethickness{1.5pt}
\put(20,0){\line(1,0){20}}
\put(20,20){\line(1,0){20}}
\put(68,10){\line(1,0){20}}
\put(68,30){\line(1,0){20}}
\put(115,0){\line(1,0){20}}
\linethickness{0.5pt}
\multiput(0,0)(7,1.4){7}%
{\line(5,1){5}}
\multiput(20,0)(7,1.4){7}%
{\line(5,1){5}}
\multiput(40,0)(7,1.4){7}%
{\line(5,1){5}}
\multiput(0,20)(7,1.4){7}%
{\line(5,1){5}}
\multiput(20,20)(7,1.4){7}%
{\line(5,1){5}}
\multiput(40,20)(7,1.4){7}%
{\line(5,1){5}}
\multiput(10,10)(7,1.4){7}%
{\line(5,1){5}}
\multiput(115,5)(7,0){3}%
{\line(1,0){5}}
\put(137,17){$H_1$}
\put(107,5){$H_2$}
\put(137,0){$H_3$}
  \end{picture}
 \end{center}
\caption{\hspace{2cm} graph $H$ \hspace{5.3cm} orbit graphs
                              \hspace{0.5cm}}
\label{fig:orbits}
  \end{figure}

We shall explain later that each orbit graph of a frame is a frame,
and each orbit graph of a median graph is $K_2$, which is
a trivial case of frames. Since
condition (ii) in Theorem~\ref{tm:main} obviously holds when $H_Q$
is $K_2$, this theorem generalizes the above result for median
metrics. On the other hand, the set of metrics $\mu$ in this
theorem does not contain some minimizable metrics since there are
frames $H$ for which (ii) is not valid. One can show that (ii) holds
when each orbit graph is either $K_2$ or $K_{2,r}$ for
$r\ge 3$, the simplest cases of frames with one orbit.
Figure~\ref{fig:orbits} illustrates the graph $H$ with three orbits,
drawn by thin, dashed and bold lines, whose orbit graphs are
$H_1\isom K_{2,3}$, $H_2\isom K_2$ and $H_3\isom K_2$.

The proof of Theorem~\ref{tm:main} will involve a number of
reductions. One of them is to show that this theorem can be derived
from Theorem~\ref{tm:mgraph} and Theorem~\ref{tm:gen_retrac} below
that claims the existence of a retraction for certain graphs.
Here a {\em retraction}\/ of a bipartite graph $K=(V(K),E(K))$
onto its subgraph $K'=(V(K'),E(K'))$ is meant to be a mapping
$\gamma:V(K)\to V(K')$
which is identical on $V(K')$ (i.e., $\gamma(v)=v$ for
all $v\in V(K')$) and brings each edge of $K$ to an edge of $K'$
(i.e., $\gamma(u)\gamma(v)\in E(K')$ for all $uv\in E(K)$).
Suppose $K$ is the Cartesian product $H_1\times\ldots\times H_k$
of graphs $H_i=(T_i,U_i)$, $i=1,\ldots,k$, i.e., $V(K)=T_1\times
\ldots\times T_k$ and nodes $(s_1,\ldots,s_k)$ and $(t_1,\ldots,t_k)$
of $K$ are adjacent if and only if there is $i\in\{1,\ldots,k\}$
such that $s_it_i\in U_i$ and $s_j=t_j$ for $j\ne i$. For a
subgraph $K'$ of $K$ and $i\in\{1,\ldots,k\}$, an $i$-{\em layer}\/
of $K'$ is a maximal subgraph of $K'$ induced by nodes
$(t_1,\ldots,t_k)$ with $t_1,\ldots,t_{i-1},t_{i+1},\ldots,t_k$
fixed.
%
%
   \begin{theorem} \label{tm:gen_retrac}
Let $K$ be the Cartesian product of frames $H_i=(T_i,U_i)$,
$i=1,\ldots,k$. Let $K'$ be an isometric subgraph of $K$ such
that $K'$ is modular and for each $i=1,\ldots,k$, some of the
$i$-layers of $K'$ is isomorphic to
$H_i$. Then there exists a retraction of $K$ onto $K'$.
   \end{theorem}
\noindent (Note that $K'$ is not an absolute retract
in general, i.e., $K'$ need not admit retraction of {\em any}\/
bipartite graph which contains $K'$ as an isometric subgraph.
Necessary and sufficient conditions on a bipartite graph to be an
absolute retract are given in \cite{BDS-87}.)
In our case, the role
of graphs $H_i$ and $K'$ in Theorem~\ref{tm:gen_retrac} will play
the graphs $H_Q$ and $H$ in Theorem~\ref{tm:main}, using the
important observation that $H$ has a canonical isometric embedding
in the Cartesian product $K$ of its orbit graphs. It turns out that
Theorem~\ref{tm:gen_retrac} can be rather easily reduced to its
special case with $k=2$; moreover, such a reduction takes place for
arbitrary modular graphs $H_1,\ldots,H_k$. To show the existence
of a retraction for this special case, with $H_1,H_2$ frames, is the
core of the whole proof of Theorem~\ref{tm:main}. Such a retraction
is just behind our ``metric uncrossing operation'', an analogue of
the cut uncrossing operation for 0-extensions of the corresponding
orbit metrics (when both $H_1,H_2$ are $K_2$, the retraction is
evident and it induces the uncrossing of two cuts, as we
explain later).

Next we deal with intractable cases. When
$\mu=d^{K_p}$, (\ref{eq:0ext}) turns into the
{\em minimum multiterminal\/ {\em(or {\em
multiway})}\/ cut problem}, which is strongly NP-hard already for
$p=3$~\cite{DJ-94}. That result
has been generalized to a larger set of path metrics.
  \begin{theorem} \label{tm:hardgr} {\rm \cite{kar-98}}
For a fixed graph $H$,
problem~(\ref{eq:0ext}) with $\mu=d^H$ is strongly NP-hard if
$H$ is non-modular or non-orientable.
  \end{theorem}

We extend this theorem as follows.
  \begin{theorem} \label{tm:nphard}
For a fixed rational-valued metric $\mu$, (\ref{eq:0ext}) is strongly
NP-hard if $\mu$ is non-modular or if the underlying graph $H(\mu)$
is non-orientable.
  \end{theorem}

The structure of this paper is as follows. Section~\ref{sec:back}
demonstrates some basic properties of modular metrics and graphs
and their orbit graphs. Section~\ref{sec:method} describes
our approach to proving Theorem~\ref{tm:main} and is aimed to explain
why this theorem reduces to Theorem~\ref{tm:gen_retrac} with $k=2$.
 \Xcomment{
The main
idea is that the result can be derived from Theorem~\ref{tm:mgraph}
if each pair $Q,Q'$ of orbits admits a retraction of the Cartesian
product of the graphs $H_Q,H_{Q'}$ onto the ``two-orbit
graph'' $H//(U-Q-Q')$. Such a
retraction can be thought of as an analog of the cut
uncrossing operation (arising in the median metric case)
for 0-extensions of the corresponding orbit metrics.
 }
The desired retraction is constructed in Section~\ref{sec:retr},
using combinatorial arguments and relying on some result concerning
the tight spans of minimizable path metrics
from~\cite{kar-98}. The construction also relies on a key lemma
proved in Section~\ref{sec:proof}. The proof of
Theorem~\ref{tm:nphard} is given in Section~\ref{sec:nphard}.
  \Xcomment{
Section~\ref{sec:concl} contains concluding
remarks and raises open questions.
}

By technical reasons, in problems~(\ref{eq:0ext}) and (\ref{eq:ext})
we will sometimes admit $\mu(st)=0$ for distinct $s,t$ and may speak
about minimizable semimetrics rather than metrics; this
does not affect the problem area in essense.
The sets of extensions and 0-extensions of a (semi)metric $\mu$ to
a set $V$ are denoted by $\exten(\mu,V)$ and
$\exten^0(\mu,V)$, respectively.
%
%
\section{Modular metrics, modular graphs, and orbits}
\label{sec:back}

By a $u$--$v$ {\em path on}\/ a set $V$ we mean any sequence
$P=(x_0,x_1,\ldots,x_k)$ of elements of $V$ with $x_0=u$ and
$x_k=v$. For a semimetric $m$
on $V$, the $m$-{\em length}\/ $m(P)$ of
$P$ is $m(x_0x_1)+\ldots+m(x_{k-1}x_k)$, and $P$ is called {\em
shortest}\/ w.r.t. $m$, or $m$-{\em shortest}, if $m(P)=m(uv)$. If
each pair $e_i=x_{i-1}x_i$ is an edge of a graph $G=(V,E)$, then
$P=(x_0,e_1,x_1,\ldots,e_k,x_k)$ is a {\em path in}\/ $G$, and we
say that $P$ is $G$-{\em shortest} if its length $|P|:=k$
is equal to $d^G(uv)$. When it is not confusing, we
abbreviate $P=x_0x_1\ldots x_k$. Given nonnegative {\em lengths}\/
$\ell(e)$ of the edges $e\in E$, we denote by $d^{G,\ell}(xy)$ the
minimum length $\ell(P)=\sum(\ell(x_{i-1}x_i): i=1,\ldots, k)$
of a path $P=x_0x_1\ldots x_k$ connecting nodes $x$ and $y$ in $G$
(the {\em path (semi)metric}\/ for $G,\ell$). From the definition
of the underlying
graph $H(\mu)$ of a metric $\mu$ it follows that $\mu=d^{H(\mu),\ell}$
for the restriction $\ell$ of $\mu$ to the edges of $H(\mu)$.

Bandelt \cite{ban-85} showed useful relations between modular
graphs and metrics. They can be stated in terms of orbits as
follows (cf.~\cite{kar-98-2}).
   \begin{myitem}
For an orbit $Q$ of a modular graph $H=(T,U)$ and nodes $u,v\in T$, if
$P$ is a shortest $u$--$v$ path and $P'$ is a $u$--$v$ path in $H$,
then $|P\cap Q|\le |P'\cap Q|$; in particular,
$|P\cap Q|=|P'\cap Q|$ if both $P,P'$ are shortest.
  \label{eq:mod1}
  \end{myitem}
   \begin{myitem}
For a modular metric $\mu$, the graph $H(\mu)$ is modular and $\mu$
is {\em orbit-invariant}, i.e., it is constant on the edges of
each orbit of $H(\mu)$.
  \label{eq:mod2}
  \end{myitem}
   \begin{myitem}
For a modular graph $H=(T,U)$ and an orbit-invariant function
$\ell:U\to\Rset_+$, the semimetric $\mu=d^{H,\ell}$ is modular,
$\mu(e)=\ell(e)$ for all $e\in U$, and every $H$-shortest path is
$\mu$-shortest; moreover, if, in addition, $\ell$ is positive, then
$H=H(\mu)$, and the metrics $d^H$ and $\mu$ have the same sets of
shortest paths.
   \label{eq:mod3}
  \end{myitem}

Note that $\mu$ need not be modular when $H(\mu)$ is modular.
(Properties~(\ref{eq:mod2}) and~(\ref{eq:mod3}) are
easily derived from~(\ref{eq:mod1}). The latter can be seen as
follows (a sketch). Let $w$ be the node of $P$ following $u$. One
may assume $P'$ is simple and all intermediate nodes $x$ of $P'$
are different from $w$. Since $P$ is shortest and $H$ is bipartite,
some node $x$ of $P'$ satisfies $d^H(wx)-1=
d^H(wy)=d^H(wz)$, where $y,z$ are the neighbours of $x$ in $P'$.
Take a median $x'$ of $y,z,w$. Then $x'y$ and $x'z$ are edges of $H$
projective to $xz$ and $xy$, respectively.
Therefore, the path $P''$ obtained from $P'$ by replacing $x$ by $x'$
obeys $|P''\cap Q|=|P'\cap Q|$, and we can apply induction since the
sum of distances from $w$ to the nodes of $P''$ is less than the
corresponding sum for $P'$, in view of $d^H(wx')=d^H(wx)-2$.)

By~(\ref{eq:mod3}), every modular graph is the underlying graph for
the class of modular metrics determined by positive
orbit-invariant functions on its edges, and all these metrics have
the same sets of shortest paths. This fact will often allow
us to work with modular graphs rather than modular metrics.

Consider a modular graph $H=(T,U)$, and let $Q_1,\ldots,Q_k$ be the
orbits of $H$. Let $\chi^S$ denote the incidence vector of a subset
$S\subseteq U$, i.e. $\chi^S(e)=1$ for $e\in S$, and 0 for $e\in U-S$.
Any modular metric $\mu$ with $H(\mu)=H$ is representable as
  \begin{equation} \label{eq:mu_sum}
  \mu=h_1\mu_1+\ldots+h_k\mu_k,
  \end{equation}
where $\mu_i=d^{H,\ell_i}$ for $\ell_i=\chi^{Q_i}$ and $h_i=\mu(e)$
for $e\in Q_i$ ($h_i$ is well-defined by (\ref{eq:mod2})).
Indeed, for any $s,t\in T$, a shortest $s$--$t$ path $P$ in
$H$ is shortest for each of $\mu,\mu_1,\ldots,\mu_k$, and $\mu_i$
coincides with $\ell_i$ on $U$, by (\ref{eq:mod3}). Therefore,
  $$
 \mu(st)=\mu(P)=h_1\ell_1(P)+\ldots+h_k\ell_k(P)=
     h_1\mu_1(st)+\ldots+h_k\mu_k(st),
  $$
as required. When all $h_i$'s are ones, (\ref{eq:mu_sum}) is
specified as
  \begin{equation} \label{eq:dH_sum}
    d^H=\mu_1+\ldots+\mu_k.
  \end{equation}

Some properties of $H$ preserve under
contraction of orbits. Let $H'=(T',U')$ be the graph $H/Q_1$. We
identify the edges in $U-Q_1$ with their images in $H'$ and denote by
$\varphi(x)$ (resp. $\varphi(P)$) the image in $H'$ of a node $x$ (resp.
a path $P$) of $H$. By~(\ref{eq:mod3}) applied to the
orbit-invariant function $\ell=\chi^{U-Q_1}$,
   \begin{myitem}
if $P$ is a shortest path of $H$, then $\varphi(P)$ is a shortest
path of $H'$.
  \label{eq:phiP}
  \end{myitem}
Therefore, if $v$ is a median of nodes $x,y,z$ in $H$, then
$\varphi(v)$ is a median of $\varphi(x),\varphi(y),\varphi(z)$
in $H'$. This implies that $H'$ is modular.
%
%
  \begin{statement} \label{st:Hprime}
$Q_2,\ldots,Q_k$ are the orbits of $H'$.
   \end{statement}
  \begin{proof}
Obviously, mates $e,e'\in U-Q_1$ of $H$ remain mates in $H'$, i.e.,
they are either opposite in a 4-circuit or parallel.
This implies that each set $Q_i$
($i>1$) is entirely included in some orbit of $H'$. To see the
reverse inclusion, consider a 4-circuit $C=(v_0,e_1,v_1,\ldots,e_4,
v_4=v_0)$ of $H'$, and let $L_j$ denote the path $(v_j,e_{j+1},
v_{j+1},e_{j+2},v_{j+2})$ for $j=0,\ldots,3$ (taking indices modulo
4). Each $L_j$ is a shortest path since $H'$ is bipartite
(as being modular).
Choose $x_0\in\varphi^{-1}(v_0)$ and $x_2\in\varphi^{-1}(v_2)$,
and let $P_0$ and $P_2$ be two $x_0$--$x_2$ paths of $H$ whose images
in $H'$ are $L_0$ and the reverse to $L_2$, respectively. Let $P$ be
a shortest $x_0$--$x_2$ path in $H$. Then $|\varphi(P)|=|L_0|=
|L_2|=2$. This together with~(\ref{eq:mod1}) (applied to $P$ and
$P'=P_0,P_2$) implies $|P\cap Q_i|=|L_0\cap Q_i|=|L_2\cap Q_i|$ for
$i=2,\ldots,k$. Similarly, $|L_1\cap Q_i|=|L_3\cap Q_i|$ for each
$i$. These equalities are possible only if each pair of mates in
$C$ belongs to the same set $Q_i$. Similar arguments
are applied to parallel edges of $H'$ (if any).
   \end{proof}

Repeatedly applying this statement to orbits of $H$, we obtain the
following.
%
%
   \begin{corollary}  \label{cor:contract}
For any $I\subseteq\{1,\ldots,k\}$, the graph $H/(\cup_{i\in I}
Q_i)$ is modular and its orbits are the sets $Q_j$ for
$j\in\{1,\ldots,k\}-I$. In particular,
each orbit graph $H_Q$ of a modular graph $H=(T,U)$ is modular and
has only one orbit, which is obtained by
identifying parallel edges in $H/(U-Q)$.
   \end{corollary}

Next we explain that each orbit graph of $H(\mu)$ is $K_2$ when $\mu$
is a median metric; this follows from properties of median graphs revealed
by Mulder and Schrijver~\cite{MS-79}. Since $\mu$ and $H(\mu)$ have
the same sets of shortest paths (by~(\ref{eq:mod2}) and (\ref{eq:mod3})),
a point $v$ is a median of points $x,y,z$ for $\mu$ if and only if
$v$ is a median of this triplet for $d^{H(\mu)}$. So $d^{H(\mu)}$ is
a median metric, which means that $H(\mu)$ is a {\em median graph}.
It is shown in~\cite{MS-79} that
   \begin{myitem}
$H=(T,U)$ is a median graph if and only if $d^H=\mu_1+\ldots+\mu_k$,
where each $\mu_i$ is the cut metric corresponding to a bi-partition
$\{A_i,T-A_i\}$ of $T$ (i.e., $\mu_i(st)=1$ if $|\{s,t\}\cap A_i|=1$,
and 0 otherwise), and the family $\Fscr=\{A_1,\ldots,A_k,
T-A_1,\ldots,T-A_k\}$ satisfies the {\em Helly property}\/
(i.e., any subfamily $\Fscr'\subseteq\Fscr$ has a nonempty
intersection provided that each two members of $\Fscr'$ meet).
  \label{eq:med_graph}
  \end{myitem}

Let $Q_i$ be the set of edges of $H$ connecting $A_i$ and $T-A_i$;
clearly $Q_1,\ldots,Q_k$ are pairwise disjoint. These
sets are precisely the orbits of $H$. Indeed,
in view of $d^H=\mu_1+\ldots+\mu_k$, a shortest path of $H$ is
$\mu_i$-shortest for each $i$. This easily implies that:
(i) the subgraphs of $H$ induced by $A_i$ and by $T-A_i$ are
connected, and (ii) $Q_i$ is a matching. (\cite{MS-79}
shows the sharper property that $H$ is median if and only if $H$ is
modular and has a cutset edge colouring.) Since $Q_i$ is
simultaneously a cut and a matching, if $e,e'$ are mates in $H$
and $e\in Q_i$, then $e'\in Q_i$. So each orbit $Q$ of $H$ is
included in some $Q_i$. Suppose $Q\ne Q_i$. Then the subgraph
$(T,U-Q)$ is connected, by (i) above, whence the semimetric
$\mu'=d^{H,\ell}$ for $\ell=\chi^Q$ is identically zero. This is
impossible because $\mu'$ coincides with $\ell$ on $U$,
by~(\ref{eq:mod3}). Thus, $Q_i$ is an orbit. Now (i) implies that
$H/(U-Q_i)$ is a tuple of parallel edges, and we conclude that
each orbit graph of $H$ is $K_2$.

As mentioned in the Introduction, our approach to solving
problem~(\ref{eq:0ext}) with a metric figured in Theorem~\ref{tm:main}
generalizes the cut uncrossing method for median metrics $\mu$.
We now briefly describe that method, referring the reader for
details to~\cite[Sec. 5]{kar-98}.

Given a median metric $\mu$ on $T$, a set $V\supseteq T$ and
a function $c:E_V\to\Zset_+$, represent $\mu$ as in~(\ref{eq:mu_sum}),
where each $\mu_i$ is the cut metric corresponding to a bi-partition
$\{A_i,T-A_i\}$ of $T$ as in~(\ref{eq:med_graph}).
For $i=1,\ldots,k$, find a bi-partition
$\{X_i,\bar X_i\}$ of $V$ such that $X_i\cap T=A_i$ and
$\sum(c(xy):x\in X_i\not\ni y)$ is minimum
(a {\em minimum cut}\/ separating $A_i$ and $T-A_i$).
Let $\Xscr=\{X_1,\ldots,X_k,\bar X_1,\ldots,\bar X_k\}$, and
let $m=h_1m_1+\ldots+h_km_k$, where $m_i$ is the cut
metric on $V$ corresponding to $\{X_i,\bar X_i\}$. Choose a pair
$Y,Z\in\Xscr$ such that $Y\cap Z\cap T=\emptyset$ but $Y,Z$ meet, and
make ``uncrossing'' by replacing $Y,Z$ in $\Xscr$ by $Y'=Y-Z$ and $Z'=Z-Y$
(taking into account that $\{Y',\bar Y'\}$ induces a minimum cut separating
$Y\cap T$ and $\bar Y\cap T$, and $\{Z',\bar Z'\}$ induces a minimum cut
separating $Z\cap T$ and $\bar Z\cap T$). Iterate until the current
family $\Xscr'$ has no such pair $Y,Z$, i.e., $Y\cap Z\cap T=\emptyset$
implies $Y\cap Z$. Using the Helly property for $\Fscr$
in~(\ref{eq:med_graph}), one can see that the corresponding metric
$m'=h_1m'_1+\ldots+h_km'_k$ is a 0-extension. Moreover,
the fact that each $m'_i$ is induced by a minimum cut implies that
$m'$ is optimal. One shows that the number of iterations does not
exceed $|T|^2 |V|$ (in fact, one can arrange a process
consisting of only $O(k^2)$ uncrossing operations).

It turns out that the Helly property for median graphs exhibited
in~(\ref{eq:med_graph}) is extended to general modular graphs. More
precisely, for a modular graph $H=(T,U)$ with orbits $Q_1,\ldots,
Q_k$, let $H_i=(T_i,U_i)$ stand for $H_{Q_i}$, and define $\pi_i=
\{A_i(t): t\in T_i\}$ to be the partition of $T$ where each member
$A_i(t)$ is the node set of the component of $(T,U-Q_i)$ whose
contraction creates the node $t$ of $H_i$. Each $A_i(t)$ is just
the corresponding maximal 0-distance set of the metric $\mu_i=
d^{H,\ell_i}$ as in~(\ref{eq:dH_sum}). We assert that
   \begin{myitem}
the family $\pi=\pi(H)$ of subsets of $T$ occurring in
$\pi_1,\ldots,\pi_k$ has the Helly property.
     \label{eq:Helly}
   \end{myitem}
Indeed, each set $A\in \pi$ is {\em convex}, i.e., for any
$x,y\in A$, each node on a shortest $x$--$y$ path $P$ of $H$
belongs to $A$. To see this, assume $A\in \pi_i$. Then $\mu_i(xy)=0$,
and therefore, $\ell_i(P)=0$ (by~(\ref{eq:mod3})). So
all nodes of $P$ are in $A$, as required.
Now the result follows from the simple fact that the family
$\bar \pi$ of convex node sets of an arbitrary modular graph
has the Helly property.
(This is shown by induction on $n$, considering a collection
$\pi'=\{A^1,\ldots,A^n\}$ of $n\ge 3$ members of $\bar\pi$
such that any $n-1$ of them meet. For $i=1,2,3$, choose an element
$x_i$ contained in all sets in $\pi'$ except possibly $A^i$. Let
$z$ be a median of $x_1,x_2,x_3$. For each $A^j\in\pi'$, at least
two of $x_1,x_2,x_3$ belong to $A^j$, hence $z\in A^j$ by the
convexity. Thus, the members of $\pi'$ have a common element.)

In conclusion of this section we show the hereditary property for
orbit graphs of frames.
   \begin{statement} \label{st:orb_frame}
Let $H=(T,U)$ be a frame, and let $Z$ be the union of some orbits
of $H$. Then $H/Z$ is a frame. In particular,
each orbit graph of $H$ is a frame.
   \end{statement}
  \begin{proof}
One can try to prove directly that the graph  $H/Z=:H'=(T',U')$
is hereditary modular and
orientable. We, however, can use Theorem~\ref{tm:mmetr} and
standard compactness arguments to show that $d^{H'}$ is
minimizable. Then $H'$ is a frame by Theorem~\ref{tm:mgraph}.

More precisely, define the semimetric $\mu$ on $T$ to be
$d^{H,\ell}$ for $\ell=\chi^{U-Z}$. Consider $V'\supseteq T'$ and
$c':E_{V'}\to\Zset_+$. We have to
show that $\tau(V',c',d^{H'})=\taustar(V',c',d^{H'})$.
Let $V=V'\cup T$ (assuming $V'\cap T=T'$) and define $c(e)=c'(e)$
for $e\in E_{V'}$, and $c(e)=0$ for
$e\in E_V-E_{V'}$. Clearly $\tau(V,c,\mu)=\tau(V',c',d^{H'})$ and
$\taustar(V,c,\mu)=\taustar(V',c',d^{H'})$. So it is suffices to
prove $\tau(V,c,\mu)=\taustar(V,c,\mu)$.

To see the latter, consider the infinite
sequence $d_1,d_2,\ldots$ of approximations for $\mu$, where
$d_i$ is $d^{H,\rho_i}$ with $\rho_i(e)=1$ for
$e\in U-Z$, and $\rho_i(e)=1/i$ for $e\in Z$. Since $H$ is
modular and $\rho_i$ is positive and orbit-invariant,
$H=H(d_i)$ for each $i$ by~(\ref{eq:mod3}).
So $d_i$ is minimizable (by Theorem~\ref{tm:mmetr}), whence
$\tau(V,c,d_i)=\taustar(V,c,d_i)$. When $i$ grows,
$\tau(V,c,d_i)$ tends to $\tau(V,c,\mu)$ (since the number of
partitions of $V$ is finite). Also $\taustar(V,c,d_i)$ tends
to $\taustar(V,c,\mu)$, because of the obvious fact that for any
$m\in \exten(\mu,V)$, there exists  $m'\in \exten(d_i,V)$
such that $|m'(e)-m(e)|\le |V|/i$ for each $e\in E_V$. Thus,
$\tau(V,c,\mu)=\taustar(V,c,\mu)$, as required.
  \end{proof}

   \Xcomment{
In conclusion of this section we point out a characterization of
modular and hereditary modular graphs that we will use later on.
  \begin{statement} \label{st:ban2} {\rm \cite{ban-88}}
{\rm (i)} A graph is hereditary modular if and only if it is bipartite
and contains no isometric $k$-circuit with $k\ge 6$.

{\rm (ii)} A modular but not hereditary modular graph
contains an isometric 6-circuit.
  \end{statement}
}
%
%
\section{Reduction to the case of two orbits, and uncrossing method}
\label{sec:method}

In this section we describe our approach to proving Theorem~\ref{tm:main}.
A majority of arguments below are applicable to general modular
metrics, and unless explicitly said otherwise, we assume that $\mu$
is an arbitrary modular metric on a set $T$.

Let $H=(T,U)$ be the underlying graph $H(\mu)$ of $\mu$ with orbits
$Q_1,\ldots,Q_k$. As before, for $i=1,\ldots,k$, $H_i=(T_i,U_i)$
stands for $H_{Q_i}$, $\ell_i$ for $\chi^{Q_i}$, $\mu_i$ for
$d^{H,\ell_i}$, and $\pi_i= \{A_i(t): t\in T_i\}$ for the corresponding
partition of $T$ defined in the previous section.
We formally identify each $t\in T_i$ with some element of
$A_i(t)$, which allows us to speak of $\mu_i$ as a 0-extension of
$d^{H_i}$ to $T$.

For the given $\mu$, consider an instance of the minimum 0-extension
problem with $V\supseteq T$ and $c:E_V\to\Zset_+$.
By~(\ref{eq:mu_sum}), any 0-extension $m$ of
$\mu$ to $V$ is representable as
  \begin{equation} \label{eq:sum_hm}
  m=h_1m_1+\ldots+h_km_k,
  \end{equation}
where each $m_i$ is the 0-extension of $\mu_i$ to $V$, defined by
   \begin{myitem}
  $m_i(xy)=\mu_i(st)$ for $x,y\in V$ and $s,t\in T$ with
$m(xs)=m(yt)=0$.
   \label{eq:def_mi}
   \end{myitem}
Then $c\cdot m=c\cdot (h_1m_1)+\ldots+c\cdot (h_km_k)$ and
$c\cdot m_i\ge\tau(V,c,\mu_i)$ for each $i$. Taking as $m$
an optimal 0-extension for $V,c,\mu$, we conclude that
  \begin{equation} \label{eq:greater}
 \tau(V,c,\mu)\ge h_1\tau(V,c,\mu_1)+\ldots +h_k\tau(V,c,\mu_k).
  \end{equation}
In particular, this is valid for $h_1=\ldots=h_k=1$ and
$\mu=d^H$. We say that $H$ is {\em orbit-additive}\/ if
  \begin{equation} \label{eq:additiv}
 \tau(V,c,d^H)=\tau(V,c,\mu_1)+\ldots +\tau(V,c,\mu_k)
  \end{equation}
holds for any $V$ and $c$. Such an $H$ has a sharper property.
%
%
  \begin{statement} \label{st:orb_add}
Let $H$ be orbit-additive. Then for any numbers $h_1,\ldots,h_k\ge 0$,
the semimetric $\mu=d^{H,\ell}$ with $\ell=h_1\ell_1+\ldots+
h_k\ell_k$ satisfies
    \begin{equation} \label{eq:h_add}
   \tau(V,c,\mu)=h_1\tau(V,c,\mu_1)+\ldots +h_k\tau(V,c,\mu_k).
    \end{equation}
Moreover, if $m$ is an optimal 0-extension for $V,c,d^H$ and
$m_1,\ldots,m_k$ are defined as in~(\ref{eq:def_mi}), then $m'=
h_1m_1+\ldots+h_km_k$ is an optimal 0-extension for $V,c,\mu$.
  \end{statement}
  \begin{proof}
Since $\tau(V,c,d^H)=c\cdot m=c\cdot m_1+\ldots+c\cdot m_k$,
(\ref{eq:additiv}) implies $c\cdot m_i=\tau(V,c,\mu_i)$ for each $i$.
Clearly $m'\in\exten^0(\mu,V)$.
Therefore, $\tau(V,c,\mu)\le c\cdot m'=h_1\tau(V,c,\mu_1)+\ldots
+h_k\tau(V,c,\mu_k)$, yielding $\tau(V,c,\mu)=c\cdot m'$ and
(\ref{eq:h_add}), in view of~(\ref{eq:greater}).
  \end{proof}

Because of~(\ref{eq:h_add}), problem~(\ref{eq:0ext}) for
a metric $\mu$ whose underlying graph $H$ is orbit-additive
becomes as easy as that for the path metrics of orbit graphs
of $H$. Indeed, to compute $\tau(V,c,\mu)$ is reduced to finding
the numbers $\tau(V,c,\mu_i)$. Moreover, once there is a subroutine
to compute $\tau(V',c',\mu)$ for arbitrary $V',c'$, we can find
an optimal 0-extension for the given $\mu,V,c$ by applying this
subroutine $O(|T| |V|)$ times (similarly to the case of minimizable
metrics $\mu$, mentioned in the Introduction).

In turn, $\tau(V,c,\mu_i)$ is equal to $\tau(V_i,c_i,d^{H_i})$,
where $V_i$ and $c_i$ arise by shrinking the sets $A_i(t)$ in the
partition $\pi_i$ of $T$ to the nodes $t\in T_i$. Formally,
$V_i=(V-T)\cup T_i$, $c_i(xy)=c(xy)$ for $x,y\in V-T$,
$c_i(xt)=c(\{x\},A_i(t))$ for $x\in V-T, t\in T_i$, and
$c_i(st)=c(A_i(s),A_i(t))$ for $s,t\in T_i$, where $c(A,B)$ denotes
$\sum(c(xy): x\in A, y\in B)$ for $A,B\subseteq V$.

In light of the above discussion, Theorem~\ref{tm:main} would follow
from Theorem~\ref{tm:mgraph} and the property that if $H$ is as in
the hypotheses of Theorem~\ref{tm:main}, then
   \begin{myitem}
   $H$ is orbit-additive.
    \label{eq:H_orb_add}
   \end{myitem}

\noindent{\bf Remark 1.} The property of being orbit-additive is
immediate in two cases of modular graphs $H$. Given $V,c$,
let $m_i$ be an optimal 0-extension for $V,c,\mu_i$,
and let $m=m_1+\ldots+m_k$. By~(\ref{eq:dH_sum}), $m\in\exten(d^H,V)$.
(i) If $H$ is a frame, then~(\ref{eq:additiv}) holds since
$\tau(V,c,d^H)=\taustar(V,c,d^H)\le c\cdot m=\tau(V,c,\mu_1)+
\ldots+\tau(V,c,\mu_k)\le \tau(V,c,d^H)$.
(ii) If $H$ is isomorphic to the Cartesian product of
$H_1,\ldots,H_k$, then $m$ is already a 0-extension of $d^H$,
yielding~(\ref{eq:additiv}); cf.~\cite{kar-98-3}.

\bigskip
We further explain that~(\ref{eq:H_orb_add}) would follow from the
existence of a retraction onto $H$ of the Cartesian product
$K=K(H)$ of the orbit graphs $H_1,\ldots,H_k$ of $H$ (see
the Introduction for needed definitions). We will use
notation $z_i$ for $i$th coordinate (component) of a point $z\in V(K)$.
Since each $H_i$ is bipartite, so is $K$. For $v\in T$,
define
   \begin{myitem}
 $\phi(v)$ to be $z\in V(K)$ such that $v\in A_i(z_i)$
for $i=1,\ldots, k$.
  \label{eq:phi}
  \end{myitem}
%
%
    \begin{statement} \label{st:embed}
For any $u,v\in T$, $d^H(uv)=d^K(\phi(u)\phi(v))$.
  \end{statement}
  \begin{proof}
Let $\phi(u)=s$ and $\phi(v)=t$. We have
$d^K(st)=d^{H_1}(s_1t_1)+\ldots +d^{H_k}(s_kt_k)$.
Consider a shortest $u$--$v$ path $P$ in $H$, and for
$i=1,\ldots,k$, let $P_i$ be the image of $P$ in $H_i$. Then
$|P|=|P_1|+\ldots+|P_k|$, and each $P_i$ is a shortest path,
by~(\ref{eq:phiP}). By~(\ref{eq:phi}),
$u\in A_i(s_i)$ and $v\in A_i(t_i)$, so
$s_i,t_i$ are the ends of $P_i$ and $|P_i|=d^{H_i}(s_it_i)$.
Therefore, $|P|=d^K(st)$.
  \end{proof}

Thus, $\phi$ induces an isometric embedding of $H$ into $K$, called
the {\em canonical}\/ embedding of $H$.
We extend $\phi$ to the edges of $H$ and, when no confusion can arise,
identify $H$ with the subgraph $\phi(H)$ of $K$. In particular, $\phi$
is injective; in other words,
   \begin{myitem}
 for $z\in V(K)$, the subset $A_1(z_1)\cap\dots\cap A_k(z_k)$
of $T$ consists of a single element (namely, $\phi^{-1}(z)$)
if $z\in\phi(T)$, and is empty otherwise.
   \label{eq:inject}
   \end{myitem}

An elementary property of a retraction of a (bipartite) graph
$G=(V,E)$ onto its subgraph $G'=(V',E')$ is that  $\gamma$ turns
every path of $G$ into a path of $G'$. This implies that
$d^G(xy)-d^{G'}(\gamma(x)\gamma(y))$ is a nonnegative even integer
for any $x,y\in V$. Therefore, $\gamma$ is {\em non-expansive}\/
(does not increase the distances) and preserves the distance
{\em parity}.
%
%
    \begin{statement} \label{st:reduction}
A modular graph $H$ is orbit-additive if there exists a retraction
of $K=K(H)$ onto $H$.
  \end{statement}
  \begin{proof}
Given $V,c$, for each $i=1,\ldots, k,$
take an optimal 0-extension $m_i$ for $V,c,\mu_i$, and form the
extension $m=m_1+\ldots+m_k$ of $d^H$ to $V$. Assuming there exists
a retraction $\gamma$ of $K$ onto $H$, we construct a
0-{\em extension}\/ $m'$ of $d^H$ to $V$ such that $m'\le m$. This
will imply~(\ref{eq:additiv}) since $\tau(V,c,\mu)\le
c\cdot m'\le c\cdot m$ and $c\cdot m=\tau(V,c,\mu_1)+\ldots+
\tau(V,c,\mu_k)$. For $z\in V(K)$, define
  \begin{eqnarray}
 X_i(z_i) &=& \{x\in V: m_i(xv)=0 \;\;\mbox{some}\;\;
  v\in A_i(z_i)\}, \;\; i=1,\ldots,k; \nonumber\\
 X_z &=& X_1(z_1)\cap\ldots \cap X_k(z_k). \label{eq:Xv}
  \end{eqnarray}

The mapping $\omega:V\to V(K)$, defined by $\omega(x)=z$ for $x\in
X_z$, isometrically embeds $(V,m)$ in $(V(K),d^K)$. Indeed, for $x\in
X_z$ and $y\in X_{z'}$, we have
  $$
 m(xy)=m_1(xy)+\ldots+m_k(xy)=d^{H_1}(z_1z'_1)+\ldots+d^{H_k}(z_kz'_k)
       =d^K(zz').
  $$
Also $\omega(v)=v$ for each $v\in T$ (cf.~(\ref{eq:phi})), i.e.,
$\omega$ is identical on the node set of the graph $H$ embedded in $K$
by $\phi$. The sets $X_z$ give a partition of $V$, and if it
happens that for each {\em nonempty} set $X_z$, the set
$A_1(z_1)\cap\ldots\cap A_k(z_k)$ is nonempty too (thus consisting
of a single node, by~(\ref{eq:inject})), then $m$ is already a
0-extension. In general, define the semimetric $m'$ on $V$ by
  $$
    m'(xy)=d^H(\gamma(\omega(x))\gamma(\omega(y)))\quad
      \mbox{for}\;\; x,y\in V.
  $$
Then $m'$ is a 0-extension of $d^H$ (corresponding to the
partition $\{\omega^{-1}\gamma^{-1}(t): t\in T\}$).
Now the fact that $\gamma$ is non-expansive while $\omega$
is isometric implies $m'\le m$, as required.
   \end{proof}

One can see that for each orbit $Q_i$, the components of the graph
$(T,Q_i)$ are just the $i$-layers of $H$ (canonically embedded in $K$
by $\phi$). Thus, condition (ii) in Theorem~\ref{tm:main} says that
each orbit graph $H_i$ is isomorphic to some of the $i$-layers of $
H$, and now summing up the above reasonings, we conclude that
Theorem~\ref{tm:main} is implied by Theorem~\ref{tm:gen_retrac}.

So it remains to prove Theorem~\ref{tm:gen_retrac}. For convenience
we denote $K'$ by $H=(T,U)$. Note that now any graph $H_i$ may have
more than one orbit, but this is not important for us. First of all
we explain that it suffices to consider the case $k=2$ (in the
reduction below we only use the fact that each $H_i$ is modular
rather than $H_i$ is a frame).

Let $1\le i<j\le k$ and $K_{ij}=H_i\times H_j$.
Define $H_{ij}=(T_{ij},U_{ij})$ to be the projection of $H$ to
$K_{ij}$, i.e., $T_{ij}=\{(z_i,z_j): z\in T\}$ and
$U_{ij}=\{(z_i,z_j)(z'_i,z'_j): zz'\in U,\; z_p=z'_p$ for $p\ne
i,j\}$. (When $H$ is as in Theorem~\ref{tm:main}, $H_{ij}$ is
isomorphic to the ``two-orbit graph'' $H//(U-Q_i-Q_j)$.)
Suppose a retraction $\gamma_{ij}$ of $K_{ij}$ onto $H_{ij}$
exists for each pair $i,j$. Define
the mapping $\psi_{ij}:V(K)\to V(K)$ by $\psi_{ij}(z)=z'$,
where $(z'_i,z'_j)=\gamma_{ij}(z_i,z_j)$ and $z'_p=z_p$ for
$p\ne i,j$. Clearly $\psi_{ij}$ is identical on $T$
and brings every edge of $K$ to an edge.
Then the desired retraction $\gamma$ of $K$ onto $H$ is
devised by successively applying transformations $\psi_{ij}$,
as follows.

At the first step, set $W_1:=V(K)$ and choose a pair
$i,j$ such that there is a point $z\in W_1$ with
$(z_i,z_j)\not\in T_{ij}$. Set $\alpha_1:=\psi_{ij}$ and reduce
$W_1$ to $W_2:=\alpha_1(W_1)$. Note that $\alpha$ decreases at
least one distance, namely, for $u=\alpha_1(z)$, we have
$\alpha_1(u)=u$, so $d^K(zu)>d^K(\alpha_1(z)\alpha_1(u))=0$. Similarly,
at each step $q$, we choose $i',j'$ with $(v_{i'},v_{j'})
\not\in T_{i'j'}$ for some $v\in W_q$, set $\alpha_q:=\psi_{i'j'}$ and
reduce $W_q$ to $W_{q+1}:=\alpha_q(W_q)$, and so on. Since each
transformation is non-expansive and brings some pair of points
of the current set $W$ to closer points, the process is finite. It
terminates when, after $N$ steps, for any $z\in W_{N+1}$, each pair
$(z_i,z_j)$ is already in $T_{ij}$. Let $\gamma=\alpha_N\alpha_{N-1}
\ldots\alpha_1$. Then $\gamma$ is identical on $T$, brings every edge
to an edge and maps $V(K)$ to $W_{N+1}$. To conclude that $\gamma$ is
a retraction of $K$ onto $H$, we have to show that $W_{N+1}=T$.
%
%
  \begin{statement} \label{st:Helly2}
Let $z$ be a point in $V(K)$ such that $(z_i,z_j)\in T_{ij}$ for all
$0\le i<j\le k$. Then $z$ is in $H$.
  \end{statement}
   \begin{proof}
For each $p=1,\ldots,k$, the set $B_p:=\{t\in T: t_p=z_p\}$ is
convex in $H$ (but not necessarily in $K$!). Indeed, if $u,v\in
B_p$ and $P$ is a shortest $u$--$v$ path in $H$, then $P$ is
shortest in $K$ (since $H$ is isometric). Therefore, $u_p=v_p=z_p$
implies $w_p=z_p$ for any node $w$ on $P$, whence $w\in B_p$.

We know that the family of convex sets of a modular graph has the
Helly property. The inclusion $(z_i,z_j)\in T_{ij}$ means that
the sets $B_i$ and $B_j$ meet. Therefore, $B_1,\ldots,B_k$ have a
common element $z'\in T$. Clearly $z'=z$.
   \end{proof}

Thus, it suffices to prove Theorem~\ref{tm:gen_retrac} for $k=2$.
The desired retraction will be constructed in the next section.

\medskip
\noindent {\bf Remark 2.} The above arguments prompt a method
to solve (\ref{eq:0ext}) with $\mu$ as in Theorem~\ref{tm:main}
in which each particular problem concerning $\mu_i$ is solved
only once (so the method looks more efficient than that
described after the proof of
Statement~\ref{st:orb_add}). More precisely, given $V,c$,
find an optimal 0-extension $m_i$ for each $i=1,\ldots,k$.
This gives the family $\Xscr$ of sets $X_i(z_i)$ as in (\ref{eq:Xv}),
and we can select, in polynomial time, the set $\Vscr$ consisting
of all points $z\in K(V)$ with $X_z\ne \emptyset$. Starting with
$\Vscr_1=\Vscr$, at each, $q$th, iteration, we examine the current set
$\Vscr_q$ to find $z\in \Vscr$ with $(z_i,z_j)\not\in T_{ij}$ for some
$i,j$. If such a $z$ exists and is chosen,
we set $\alpha_q:=\psi_{ij}$, reduce $\Vscr_q$ to $\Vscr_{q+1}:=
\alpha_q(\Vscr_q)$ (which changes $\Xscr$)
and continue the process. Otherwise $\Vscr_q=T$,
by Statement~\ref{st:Helly2}, and the partition $\{Y_t\>:\;t\in T\}$
of $V$ into the corresponding 0-distance sets induces an optimal
0-extension for $V,c,d^H$ (and therefore, for $V,c,\mu$, by
Statement~\ref{st:orb_add}), where $Y_t$ is the union
of sets $X_z$ for $z\in\Vscr$ such that
$\alpha_{q-1}\ldots\alpha_1(z)=t$. Since each transformation
moves some point of the current set $\Vscr$
closer to $T$, the number of iterations is $O(|T|^2|V|)$.

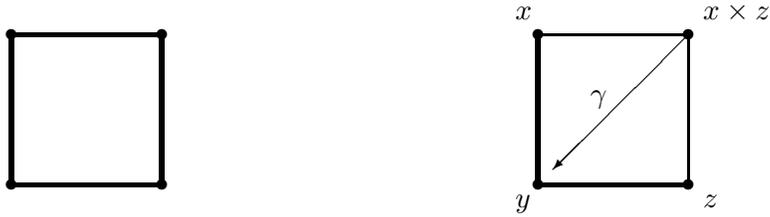
\begin{figure}[tb]
  \unitlength=1mm
 \begin{center}
  \begin{picture}(125,25)(-20,5)
\put(5,5){\circle*{1.5}}  
\put(25,5){\circle*{1.5}}
\put(5,25){\circle*{1.5}}
\put(25,25){\circle*{1.5}}
\linethickness{1.5pt}
\put(5,5){\line(0,1){20}}
\put(5,5){\line(1,0){20}}
\put(5,25){\line(1,0){20}}
\put(25,5){\line(0,1){20}}
\put(75,5){\circle*{1.5}}  
\put(95,5){\circle*{1.5}}
\put(75,25){\circle*{1.5}}
\put(95,25){\circle*{1.5}}
\put(75,5){\line(0,1){20}}
\put(75,5){\line(1,0){20}}
\linethickness{0.5pt}
\put(75,25){\line(1,0){20}}
\put(95,5){\line(0,1){20}}
\put(95,25){\vector(-1,-1){18}}
\put(72,2){$y$}
\put(97,2){$z$}
\put(72,27){$x$}
\put(97,27){$x\times z$}
\put(82,16){$\gamma$}
  \end{picture}
 \end{center}
\caption{\hspace{1cm} (a) $H_{ij}\simeq K_2\times K_2$ \hspace{4cm}
  (b) $H_{ij}\simeq P$ \hspace{2cm} } \label{fig:Hij}
  \end{figure}

\medskip
\noindent {\bf Remark 3.}
The above transformation of $\Xscr$ induced by the retraction
$\gamma_{ij}$ can be thought of as an analogue of the cut
uncrossing operation for median metrics
(reviewed in Section~\ref{sec:back}), thus justifying the term
``uncrossing'' used in a more general context in this paper. Recall
that each orbit graph of a median graph $H$ is $K_2$, and therefore,
each ``two-orbit graph'' $H_{ij}$ is isomorphic
either to $K_2\times K_2$ or to the path $P=xyz$ of length two,
as drawn in Fig.~\ref{fig:Hij}.
When $H_{ij}\simeq P$, the (unique) retraction $\gamma=\gamma_{ij}$
brings the point $(x,z)$ of $H_i\times H_j$ not in $H_{ij}$ to $y$.
This retraction is just behind the uncrossing operation on
the corresponding cuts in that method.

%
%
\section{Retraction} \label{sec:retr}

In this and next sections we prove Theorem~\ref{tm:gen_retrac}
with $k=2$, using notation, conventions and results from
Sections~\ref{sec:back} and \ref{sec:method}.
One may assume $K\ne H$. We will essentially use the condition in
the theorem that $H$ includes a subgraph (``row-layer'') of the
form $H_1\times s_2$ and a subgraph (``column-layer'')
of the form $s_1\times H_2$ for some $s_1\in T_1$ and $s_2\in T_2$,
i.e.,
   \begin{myitem}
   for any $u\in T_1$ and $v\in T_2$, $(u,s_2)\in T$ and
 $(s_1,v)\in T$.
   \label{eq:origin}
   \end{myitem}
We fix such $s_1,s_2$ and call the node $s=(s_1,s_2)$ the
{\em origin}\/ of $K$.

In the proof below we everywhere admit that $H_1,H_2$ are
arbitrary modular graphs until~(\ref{eq:tight_span})
where the assumption that $H_1,H_2$ are frames is essential.
We abbreviate $d^K,d^{H_1},d^{H_2}$ as $d,d_1,d_2$,
respectively. The {\em interval}\/ $\{v\in V(K): d(xv)+d(vy)=d(xy)\}$
of nodes (points) $x,y$ of $K$ is denoted by $I(x,y)=I(y,x)$. We
denote by $J(x)$ and $r(x)$ the interval $I(x,s)$ and the distance
$d(xs)$, called the {\em principal}\/ interval and the {\em rank}\/
of $x$, respectively. $M(x,y,z)$ denotes the set of medians of
points $x,y,z\in V(K)$. For $i=1,2$, $I_i(x_i,y_i)$,
$J_i(x_i)$, $r_i(x_i)$, and $M_i(x_i,y_i,z_i)$ stand for the
corresponding objects concerning the graph $H_i$. Then
$I(x,y)=I_1(x_1,y_1)\times I_2(x_2,y_2)$, $J(x)=J_1(x_1)\times
J_2(x_2)$, $r(x)=r_1(x_1)+r_2(x_2)$, and $M(x,y,z)=M_1(x_1,y_1,z_1)
\times M_2(x_2,y_2,z_2)$ (as being immediate consequences from
the equality $d(uv)=d_1(u_1v_1)+d_2(u_2v_2)$ for any $u,v\in V(K)$).
The latter correspondence between medians in $K,H_1,H_2$ implies the
following elementary property, which will be often used later on:
    \begin{myitem}
  for $x,y,z\in V(K)$ and $i\in\{1,2\}$, if $z_i\in
I_i(x_i,y_i)$, then $z_i=v_i$ for each median $v\in M(x,y,z)$;
in particular, $x_i=z_i$ implies $v_i=x_i$.
    \label{eq:med_int}
    \end{myitem}

The modularity of $H$ implies that
    \begin{myitem}
  for each $u\in T_1$, the set $Z(u):=\{v\in T_2: (u,v)\in T\}$ is
convex in $H_2$, and similarly for each $v\in T_2$, the set
$\{u\in T_1: (u,v)\in T\}$ is convex in $H_1$
    \label{eq:conv_Z}
    \end{myitem}
(cf. the proof of Statement~\ref{st:Helly2}).
Indeed, for $v,w\in Z(u)$ and $v'\in I_2(v,w)$, consider
the nodes $x=(u,v)$, $y=(u,w)$ and $z=(s_1,v')$ of $H$ (where $z$ is
in $T$ by~(\ref{eq:origin})). These nodes have a median
$q$ in $H$. Then $q_1=u$ and $q_2=v'$ (cf.~(\ref{eq:med_int})).
Hence, $v'\in Z(u)$. It follows from~(\ref{eq:conv_Z})
that
    \begin{equation} \label{eq:local_conv}
   J(t)\subseteq T\qquad \mbox{for all}\;\; t\in T.
    \end{equation}
(However, the whole set $T$ is not convex in $K$ unless $H=K$.)

The mapping (retraction) $\gamma$ that we wish to construct will be
some kind of reflection of points in $V(K)-T$ with respect to
their closest points in $H$. Consider a point $x\in V(K)$. Define
the {\em excess}\/ $\Delta^x$ to be the distance $d(x,T)$ from $x$ to $T$,
i.e., $\Delta^x=\min\{d(xt): t\in T\}$, and define $N(x)$ to be the
set of points $t\in T$ with $d(xt)=\Delta^x$. In particular,
$\Delta^x\le r_i(x_i)$ for $i=1,2$ (since $(x_1,s_2),(s_1,x_2)\in T$),
and $\Delta^x=0$ if and only if $x\in T$.
%
%
  \begin{statement} \label{st:N_J}
 $N(x)\subseteq J(x)$.
  \end{statement}
 \begin{proof}
Let $t\in N(x)$. The points $x'=(x_1,s_2)$, $x''=(s_1,x_2)$ and $t$
are in $T$, so they have a median $q$ in $T$ as well. Then $q_1\in
M_1(x_1,s_1,t_1)$ and $q_2\in M_2(s_2,x_2,t_2)$. Therefore, $q_1$
belongs to both $J_1(x_1)$ and $I_1(x_1,t_1)$, and $q_2$ belongs to
both $J_2(x_2)$ and $I_2(x_2,t_2)$, which means that $q\in J(x)$
and $q\in I(x,t)$. Now $d(xq)\ge\Delta^x=d(xt)$ implies $q=t$.
  \end{proof}

By this statement, the rank $r(t)$ is equal to the same number
$r(x)-\Delta^x$ for all $t\in N(x)$. Note that for any $x,y\in V(K)$,
$|\Delta^x-\Delta^y|=|d(x,T)-d(y,T)|\le d(xy)$. Therefore,
   \begin{equation} \label{eq:Delta}
  |\Delta^x-\Delta^y|\le 1 \qquad\mbox{for each edge}\;\;\;
        xy\in E(K).
    \end{equation}

We partition $E(K)$ into the sets $E_1=\{xy: x_2=y_2\}$
and $E_2=\{xy: x_1=y_1\}$, and for $i=1,2$, define
  \begin{equation} \label{eq:E1E2}
    \Eie=\{xy\in E_i: \Delta^x=\Delta^y\}\quad \mbox{and}\quad
       \Eine=E_i-\Eie.
   \end{equation}

The desired retraction is devised by use of certain
0-extensions of metrics $d_1$ and $d_2$. First we introduce the
auxiliary graphs $G_1=(\Vscr_1,\Escr_1)$ and $G_2=
(\Vscr_2,\Escr_2)$, as follows.
For $i=1,2$, let $\Ascr_i$ be the set of pairs $tt_i=\{t,t_i\}$ for
$t\in T$, and $\Bscr_i$ the set of pairs $xs_i=\{x,s_i\}$ for
$x\in V(K)$. Then $G_i$ is the (disjoint) union of the graphs $H_i$
and $K$ to which the pairs from $\Ascr_i\cup\Bscr_i$ are added as edges,
i.e.,
   $$
  \Vscr_i=T_i\cup V(K)\quad \mbox{and}\quad
    \Escr_i=U_i\cup E(K)\cup\Ascr_i\cup\Bscr_i.
   $$
The edges $e$ of $G$ are endowed with the lengths $\delta_i(e)$
defined by
     \begin{eqnarray}
  \delta_i(e) &=& 1 \quad \mbox{for} \quad e\in U_i\cup \Eie
                    \cup E^{\ne}_{3-i},  \label{eq:delta_i}\\
     &=& 0 \quad \mbox{for}\quad e\in\Eine\cup E^=_{3-i} \cup
                  \Ascr_i,  \nonumber  \\
     &=& r_i(x_i)-\Delta^x \quad \mbox{for}\quad
                       e=xs_i\in \Bscr_i. \nonumber
  \end{eqnarray}

We say that a semimetric $m$ on a set $V$ is {\em cyclically even}\/
if $m(xy)+m(yz)+m(zx)$ is an even integer for all $x,y,z\in V$
(equivalently: the $m$-length of any cycle on $V$ is even).
All values of such an $m$ are integers since $m(xy)+m(yx)+m(xx)
=2m(xy)\in 2\Zset$.

%
%
    \begin{lemma} \label{lm:m_i}
For $i=1,2$, define
$m_i=d^{G_i,\delta_i}$. Then: {\rm (i)}\/ $m_i$ is an extension of
$d_i$ to $\Vscr_i$, and {\rm (ii)}\/ $m_i$ is cyclically even and
coincides with $\delta_i$ on $\Escr_i$.
    \end{lemma}
This lemma (the keystone in our arguments) will be proved later,
and now we explain how it help us to construct the desired mapping
$\gamma$. We apply some results from~\cite{kar-98-3} and
\cite{kar-98}.

More precisely, for a metric $\mu'$ on a set $T'$, an extension
$m'$ of $\mu'$ to $V\subseteq T'$ is called {\em tight}\/ if there
exists no $m''\in \Escr(\mu',V)-\{m'\}$ such that $m''\le m'$;
equivalently: $m'$ has no {\em loose}\/ pair $x,y$, i.e, for any
$x,y\in V$, the path $(u,x,y,v)$ on $V$ is
$m'$-shortest for at least one pair $u,v\in T'$.
It is shown in~\cite[Sec.5]{kar-98-3}
that for any cyclically even metric $\mu'$,
   \begin{myitem}
  if $m\in\exten(\mu',V)$ is cyclically even, then there exists
$m'\in\exten(\mu',V)$ such that $m'$ is cyclically even and tight,
$m'(e)\le m(e)$ for all $e\in E_V$, and $m'(e)=m(e)$ whenever
$m(e)\le 1$.
    \label{eq:c_e_tight}
   \end{myitem}
(Such an $m'$ is constructed by the following process. If there is
no loose pair $x,y\in V$ with $m(xy)\ge 2$, then one easily shows that
there is no loose pair at all, i.e., $m$ is already tight.
Otherwise choose such a pair $x,y$, and let $m':=d^{K_V,\ell}$,
where $\ell(xy):=m(xy)-2$ and $\ell(e):=m(e)$ for $e\in E_V-\{xy\}$.
Then $m'$ is a cyclically even extension of $\mu'$. Update $m:=m'$
and iterate.)

Next, the proof of the ``if'' part of Theorem~\ref{tm:mgraph} in
\cite{kar-98} relies of an explicit construction of the so-called
tight span of a frame, which in turn is based on the following
result (Claim 5 in Section 4 there):
   \begin{myitem}
 if $H'=(T',U')$ is a frame and $m$ is a tight extension of
$d^{H'}$ to $V\supseteq T'$, then each point $x\in V$ satisfies
at least one of the following:
    \begin{itemize}
   \item[(i)] $m(xt)=0$ for some node $t\in T'$;
   \item[(ii)] $m(ux)+m(xv)=1$ for some edge $uv\in U'$;
   \item[(iii)] $m(v_0x)+m(xv_2)=m(v_1x)+m(xv_3)=2$ for some
4-circuit $C=v_0v_1v_2v_3v_0$ of $H'$.
     \end{itemize}
  \label{eq:tight_span}
    \end{myitem}

Using~(\ref{eq:c_e_tight}) and (\ref{eq:tight_span}), we argue as
follows. For $i=1,2$, let $m_i$ be as in Lemma~\ref{lm:m_i}, and let
$m'_i\le m_i$ be a cyclically even tight extension of $d_i$ as
in~(\ref{eq:c_e_tight}). Then
   \begin{equation}  \label{eq:m_and_delta}
  m'_i(e)=\delta_i(e)\;\; \mbox{for}\;\; e\in \Escr_i-\Bscr_i,
  \quad\mbox{and}\quad m'_i(e)\le\delta_i(e)\;\; \mbox{for}\;\;
          e\in\Bscr_i.
   \end{equation}
Moreover, in view of (\ref{eq:tight_span}), for each $x\in\Vscr_i$,
there exists $t\in T_i$ with $m'_i(tx)=0$. This is immediate in
cases (i) and (ii) of~(\ref{eq:tight_span}). And if we are in case (iii)
(with $m=m'_i$) and if $C=v_0v_1v_2v_3v_0$ is the corresponding 4-circuit
for $x$, then $\alpha_j:=m'_i(v_jx)>0$ for $j=0,1,2,3$ would imply
$\alpha_j=1$ for each $j$. Then $m'_i(v_0v_1)+\alpha_0+\alpha_1=
1+1+1=3$, contrary to the fact that $m'_i$ is cyclically even. Thus,
$m'_i$ is a 0-extension of $d_i$ to $\Vscr_i$.

Now for $x\in V(K)$, define $\gamma(x)$ to be the point
$(\gamma_1(x),\gamma_2(x))$, where $\gamma_i(x)$ is the node $v\in T_i$
with $m'_i(xv)=0$.
%
%
   \begin{statement} \label{st:the_gamma}
  $\gamma$ is the retraction of $K$ onto $H$.
   \end{statement}
  \begin{proof}
 For each $t\in T$, $m'_i(tt_i)=0$ (since $\delta_i$ is zero on
$\Ascr_i$, by~(\ref{eq:delta_i})), so $\gamma$ is identical on $T$.

To see $\gamma(V(K))\subseteq T$, consider $x\in V(K)$, and let
$x'=\gamma(x)$ and $t\in N(x)$. Let $P=z^0z^1\ldots z^k$ ($k=\Delta^x$)
be a shortest $t$--$x$ path in $K$. Then for $j=0,\ldots,k-1$,
one has $t\in N(z^j)$ and $\Delta_j:=\Delta^{z^j}=j$, whence
$\Delta_{j+1}-\Delta_j=1$ and $z^jz^{j+1}\in \Eonene\cup\Etwone$,
cf.~(\ref{eq:E1E2}). This implies $\delta_1(P)=d_2(t_2x_2)$ and
$\delta_2=d_1(t_1x_1)$, by the definition of $\delta_i$ on $E(K)$.
Therefore,
   \begin{equation} \label{eq:first}
 d_1(x'_1t_1)=m'_1(xt)\le\delta_1(P)=d_2(t_2x_2)=\Delta^x-d_1(t_1x_1).
    \end{equation}
Since $\delta_1(s_1x)=r_1(x_1)-\Delta^x$ (by~(\ref{eq:delta_i})) and
$r_1(x_1)=r_1(t_1)+d_1(t_1x_1)$ (by Statement~\ref{st:N_J}),
   \begin{equation} \label{eq:second}
 d_1(s_1x'_1)=m'_1(s_1x)\le\delta_1(s_1x)=r_1(x_1)-\Delta^x
  =r_1(t_1)+d_1(t_1x_1)-\Delta^x.
   \end{equation}

Comparing~(\ref{eq:first}) and (\ref{eq:second}), we obtain
$d_1(s_1x'_1)+d_1(x'_1t_1)\le r_1(t_1)$, whence $x'_1\in J_1(t_1)$.
Similarly, $x'_2\in J_2(t_2)$. So $x'\in J(t)$, yielding $x'\in T$,
by~(\ref{eq:local_conv}).

Finally, consider an edge $e=xy\in E(K)$, and let $x'=\gamma(x)$ and
$y'=\gamma(y)$. We have $\delta_1(e)+\delta_2(e)=1$,
by~(\ref{eq:delta_i}). Also $m'_i(e)=\delta_i(e)$, $i=1,2$,
by~(\ref{eq:m_and_delta}). Hence,
  $$
   d(x'y')=d_1(x'_1y'_1)+d_2(x'_2y'_2)=m'_1(e)+m'_2(e)=
     \delta_1(e)+\delta_2(e)=1,
  $$
i.e., $x'y'$ is an edge of $K$, as required.
  \end{proof}

It remains to prove Lemma~\ref{lm:m_i}.

%
%
\section{Proof of Lemma~\ref{lm:m_i} } \label{sec:proof}

We may prove this lemma for $i=1$. First
we explain that $\delta_1$ is cyclically even, i.e., the
$\delta_1$-length of any cycle in $G_1$ is even.

For any 4-circuit $C=x^0x^1x^2x^3x^0$ in $K$, an edge of $C$ belongs
to $E_1$ if and only if the opposite edge does. Also, letting
$\eta_j:=\Delta^{x^{j+1}}-\Delta^{x^j}$, the numbers
$\eta_0,\eta_2$ have the same parity if and only if $\eta_1,\eta_3$
do so. From these properties and the definition of $\delta_i$
one can deduce that the $\delta_1$-length of
$C$ is even. Then $\delta_1$ is cyclically even within $K$, because
$K$ is modular and, therefore, the 4-circuits form a basis in the
space of cycles of $K$ over $\Zset_2$. (Indeed, any cycle of length
$q\ge 6$ in a modular graph can be represented as the modulo
two sum of three cycles with length less than $q$ each.) Next,
using the fact that $\delta_1$ takes value one on
$U_1\cup(E_1\cap U)$ and zero on $(E_2\cap U)\cup\Ascr_1$, one can
see that the $\delta_1$-length of any cycle with all edges in
$U_1\cup U\cup\Ascr_1$ is even. Finally, for an edge
$e=xs_1\in\Bscr_1$, choose $t\in N(x)$ and a shortest $t$--$x$ path
$L$ in $K$. Then $\delta_1(L)=d_2(t_2x_2)$. Concatenating $L$
with the edge $e$,
the edge $t_1t$ in $\Ascr_1$ and a shortest $s_1$--$t_1$ path $R$
in $H_1$, we obtain a cycle whose $\delta_1$-length is equal to
  $$
 \delta_1(L)+\delta_1(e)+\delta_1(R)+\delta_1(t_1t)=
  d_2(t_2x_2)+(r_1(x_1)-\Delta^x)+r_1(t_1)+0=2r_1(t_1).
  $$

Summing up the above observations, one can conclude that $\delta_1$ is
cyclically even within the entire set $\Escr_1$. Then
$m_1$ is cyclically even as well.

Next we prove that $m_1$ is an extension of $d_1$. The main part of
this proof is to show the following property:
    \begin{myitem}
for any path $P=x^0x^1\ldots x^k$ in $K$ with $x^0\in T$, there
exists a path $L=z^0z^1\ldots z^\alpha$ with $z^0=x^0$ and
$z^\alpha=x^k$ and a number $0\le \beta\le \alpha$ such that
$z^0,\ldots,z^\beta\in T$, that $r(z^\beta)<r(z^{\beta+1})<
\ldots<r(z^\alpha)$, and that $\delta_1(L)\le \delta_1(P)$.
    \label{eq:path_P}
   \end{myitem}

The proof of (\ref{eq:path_P}) includes Claims 1--3 below.
Recall that any edge $xy\in E(K)$
satisfies $|r(x)-r(y)|=1$ (since $K$ is bipartite), and if
$x\in T$ and $r(x)>r(y)$, then $y\in T$ (by~(\ref{eq:local_conv})).
In particular, $L$ as in~(\ref{eq:path_P}) entirely lies in $H$
if $x^k\in T$. To show~(\ref{eq:path_P}), it suffices
to consider the case when $P$ is simple, $k\ge 2$, and all
intermediate nodes of $P$ are not in $T$ (for if $x^i\in T$ for
some $0<i<k$, we can split $P$ into two paths $P'=x^0\ldots x^i$
and $P''=x^i\ldots x^k$ and prove~(\ref{eq:path_P}) for each of
$P',P''$ independently). For $i=0,\ldots,k$, let $r(i):=r(x^i)$.
An intermediate node $x^i$ of $P$ is
called a {\em peak}\/ if $r(i)>r(i-1)=r(i+1)$. The set of peaks is
denoted by $F=F(P)$. We prove~(\ref{eq:path_P}) by induction on
   $$
  \omega(P)=\sum(4^{r(i)}: x^i\in F(P)).
   $$

If $F=\emptyset$, then $r(0)<r(1)<\ldots<r(k)$ (as $r(0)>r(1)$
would imply $x^1\in T$), i.e., $P$ is just the desired path $L$.
So assume $F\ne\emptyset$. Let $x^p$ be the first peak in $P$, and
let $x,y,z$ stand for $x^{p-1},x^p,x^{p+1}$, respectively. Choose a
median $y'$ for $x,z,s$ in $K$. Since $r(x)=r(z)$ and $d(xz)=2$,
both $xy',y'z$ are edges of $K$ and $r(y')<r(x)<r(y)$.
Replace $y$ by $y'$ in $P$, forming the path $P'=x^0\ldots x^{p-1}
y'x^{p+1}\ldots x^k$; we say that $P'$ is obtained by {\em
cutting off}\/ the peak $y$. Since $4^{r(p)}>2\cdot 4^{r(p)-1}=
4^{r(p-1)}+4^{r(p+1)}$, we have $\omega(P')<\omega(P)$. Also
$\delta_1(P)-\delta_1(P')$ is equal to
   $$
  \rho:=\rho(x,y,z,y'):= \delta_1(xy)+\delta_1(yz)-
       \delta_1(xy')-\delta_1(y'z).
   $$
Therefore, if $\rho\ge 0$ occurs, we can immediately apply
induction. Let $\bar\Delta:=\Delta^y$.

\bigskip
\noindent {\bf Claim 1} {\em A median $y'$ for $x,z,s$ can be chosen
so that $\rho(x,y,z,y')<0$ is possible only if both edges
$e=xy,e'=yz$ are in $E_2$, $\Delta^{x}=\Delta^{z}=\bar\Delta$,
and $\Delta^{y'}=\bar\Delta-1$. }

\bigskip
 \begin{proof} Since the $\delta_1$-length of the 4-circuit
$C=xyzy'x$ is even, $\rho<0$ implies
    \begin{equation} \label{eq:delta_on_C}
 \delta_1(e)=\delta_1(e')=0 \qquad\mbox{and} \qquad
  \delta_1(xy')=\delta_1(y'z)=1.
   \end{equation}
This is impossible when $e\in E_1$ and $e'\in E_2$ (or $e\in E_2$
and $e'\in E_1$). Indeed, in this case we would
have $\Delta^x=\bar\Delta-1$ and $\Delta^z=\bar\Delta$,
by~(\ref{eq:delta_i}). Then $xy'\in E_2$ and $\delta_1(xy')=1$ imply
$\Delta^{y'}=\Delta^x-1=\bar\Delta-2$, while $y'z\in E_1$ and
$\delta_1(y'z)=1$ imply $\Delta^{y'}=\Delta^z=\bar\Delta$; a
contradiction.

If $e,e'\in E_2$, then $xy',y'z\in E_2$. So~(\ref{eq:delta_on_C})
yields $\Delta^x=\Delta^z=\bar\Delta$ and
$\Delta^{y'}=\Delta^x-1=\bar\Delta-1$, as required.

Now, suppose $e,e'\in E_1$ and $\delta_1(e)=\delta_1(e')=0$.
Choose $u\in N(x)$ and $v\in N(z)$. We have $\Delta^x=\Delta^z=\bar
\Delta-1$, whence $u,v\in N(y)$. Choose in $T$ a median $q$ for
$u,v,(y_1,s_2)$ and a median $w$ for $u,v,(s_1,y_2)$. We assert that
$q,w\in N(y)$. Indeed,
   $$
q_1\in M_1(u_1,v_1,y_1),\quad w_1\in M_1(u_1,v_1,s_1),\quad
q_2\in M_1(u_2,v_2,s_2),\quad w_2\in M_1(u_2,v_2,y_2).
   $$
In particular, $q_1,w_1\in I_1(u_1,v_1)$. Also
$u_1,v_1\in I_1(q_1,w_1)$ (in view of $u_1,v_1\in I_1(y_1,s_1)$,
by Statement~\ref{st:N_J}). These relations imply $d_1(u_1q_1)=
d_1(v_1w_1):=a$. Similarly, $d_2(u_2q_2)=d_2(v_2w_2):=a'$. Then
$d(yu)=\bar \Delta\le d(yq)=d(yu)-a+a'$ and $d(yv)\le d(yw)=
d(yv)+a-a'$. This is possible only if $a=a'$, yielding
$d(yq)=d(yw)=\bar\Delta$, as required.

Assume $y'$ is chosen to be a median for $x,z,w$. Then $y'$ is a median for
$x,z,s$ as well, taking into account that $x_2=z_2$ and the paths
$(x_1,u_1,w_1,s_1)$ and $(z_1,v_1,w_1,s_1)$ on $T_1$ are
$d_1$-shortest. Now $d(y'w)=d(xw)-1$ implies $\Delta^{y'}<\Delta^x$.
Hence, $\delta_1(xy')=\delta_1(y'z)=0$ and $\rho=0$.
 \end{proof}

Arguing as in the above proof, one can see that for any
$x'\in V(K)$, there are elements $t,t'\in N(x')$ such that
$r_1(t_1)\le
r_1(t'_1)$ (and $r_2(t_2)\ge r_2(t')$) and $N(x')\subseteq I(t,t')$.
We denote $t$ by $t(x')$ and refer to it as the {\em minimal}\/
element of $N(x')$ (with respect to the rank in $H_1$).

\bigskip
\noindent {\bf Remark 4.} For $i=1,2$ and $f,g\in N(x')$, denote
$f_i\prec_i g_i$ if $f_i\in J_i(g_i)$. Then $\prec_i$ is the
partial order on $N_i=\{w_i: w\in N(x')\}$ with unique
minimal and maximal elements. Moreover, the correspondence
$w_1\to w_2$
establishes the isomorphism between $(N_1,\prec_1)$ and
$(N_2,\prec_2^{-1})$ (where $\prec^{-1}$ is the reverse to $\prec$).
One can show that if none of $H_1,H_2$ containes $\kott$ as an
induced subgraph (see Fig.~\ref{fig:orient}b), then $(N_i,\prec_i)$
is a {\em modular lattice}, i.e., (i) any $u,v\in N_i$ have
unique lower and upper bounds, denoted by $u\wedge v$ and $u\vee v$,
respectively; (ii) for each $u\in N_i$, all maximal chains to $u$
from the minimal element have the same length $\rho(u)$, and (iii)
each pair $u,v$ satisfies the modular equality $\rho(u)+\rho(v)=
\rho(u\wedge v)+\rho(u\vee v)$.
We, however, do not need these properties in further arguments.

\bigskip
In light of Claim 1, we may assume that $\rho<0$ and $e,e'\in E_2$.
Consider the minimal element $t(y)=(t_1(y),t_2(y))$ in $N(y)$.
Suppose $t_1(y)\ne y_1$. Then there is a node $w$ of $K$ adjacent to
$y$ such that $w_1\in I_1(y_1,t_1(y))$ and $w_2=y_2$. We have
$yw\in E_1$, $r(w)=r(y)-1$ and $t(y)\in N(w)$. Then
$\Delta^w<\bar\Delta$ and $\delta_1(yw)=0$. Transform $P$
into the (non-simple) path $P'=x^0\ldots x^{p-2}xywyzx^{p+2}
\ldots x^k$ and then cut off both copies of $y$ (which are
peaks of $P'$). This results in a path $P''$ of the form $x^0\ldots
x^{p-2}xy'wy''zx^{p+2}\ldots x^k$; clearly $x,w,z$ are peaks of $P''$.
Since $yw\in E_1$, $y'$ and $y''$ can be chosen so that
$\rho(x,y,w,y')\ge 0$ and $\rho(w,y,z,y'')\ge 0$, by Claim 1. Therefore,
$\delta_1(P'')\le \delta_1(P')=\delta_1(P)$. Also $4^{r(y)}>
3\cdot 4^{r(y)-1}=4^{r(x)}+4^{r(w)}+4^{r(z)}$, yielding
$\omega(P'')<\omega(P)$. So we can apply induction.

It remains to consider the case when $t_1(y)=y_1$. Then $t(y)$ is
the unique element of $N(y)$. We will use the following property.

\bigskip
\noindent {\bf Claim 2.} {\em Let $\bar x\bar y\in E_2$ satisfy
$r(\bar x)<r(\bar y)$, let $N(\bar y)$ consist of a single element
$u$, and let $u_1=\bar y_1$. Then $N(\bar x)$ consists of a single
element $v$, and $v_1=\bar y_1$. Moreover, $u=v$ if $\Delta^{\bar x}
<\Delta^{\bar y}$, and $u$ and $v$ are adjacent if $\Delta^{\bar x}
=\Delta^{\bar y}$.  }

\medskip
  \begin{proof}
If $\Delta^{\bar x}<\Delta^{\bar y}$, then $N(\bar x)\subseteq
N(\bar y)$, whence $N(\bar x)=\{u\}$. So assume $\Delta^{\bar x}=
\Delta^{\bar y}$, and let $v\in N(\bar x)$. Choose $q\in M(u,v,
(\bar y_1,s_2))\cap T$ and $w\in M(u,v,(s_1,\bar y_2))\cap T$. We
have $q_2,w_2\in I_2(u_2,v_2)$ and $u_2\in I_2(q_2,w_2)$ (in view
of $u_2\in I_2(\bar y_2,s_2)$). Note that the path $(\bar y_2,\bar x_2,
v_2,s_2)$ on $T_2$ is $d_2$-shortest (since $r(\bar x)<r(\bar y)$
and $\bar x_1=\bar y_1$ imply $\bar x_2\in I_2(\bar y_2,s_2)$). This
yields $v_2\in I_2(q_2,w_2)$, and we can conclude that $d_2(u_2w_2)
=d_2(v_2q_2)=:a'$.

Next, $q_1\in M_1(u_1,v_1,\bar y_1)$ and $u_1=\bar y_1$ imply
$q_1=\bar y_1$, while $w_1\in M_1(u_1,v_1,s_1)$,
$v_1\in I_1(\bar x_1,s_1)$ and $\bar x_1=\bar y_1=\bar u_1$ imply
$w_1=v_1$. Let $a:=d_1(\bar y_1v_1)$. Then $d(\bar xv)\le d(\bar xq)
=d(\bar xv)-a+a'$ and $d(\bar yu)\le d(\bar yw)=d(\bar yu)+a-a'$,
whence $a=a'$, $q\in N(\bar x)$ and $w\in N(\bar y)$. Since
$|N(\bar y)|=1$, we have $w=u$. This implies $a=0$ and $q=v$, yielding
$v_1=q_1=\bar y_1$. So $v_1=\bar y_1$, regardless of the choice of
$v$ in $N(\bar x)$. This is possible only if $N(\bar x)$ consists of
a single element (for if $v,v'\in N(\bar x)$ and $v\ne v'$, then
a median $f$ for $v,v',(s_1,\bar x_2)$ in $T$ satisfies $f_1=\bar y_1$
and $d_2(\bar x_2f_2)<d_2(\bar x_2v_2)$, whence
$d(\bar xf)<\Delta^{\bar x}$).

Finally, to see that $u_2$ and $v_2$ are adjacent, take in $T$ a median
$h$ for $u,v,(s_1,\bar x_2)$. Then $d(\bar xh)\ge d(\bar xv)$,
$h_2\in I_2(\bar x_2,v_2)$ and $h_1=\bar y_1$, implying $h=v$. So
$v_2\in I_2(\bar x_2,u_2)$. Also $d_2(\bar y_2u_2)=\Delta^{\bar y}=
\Delta^{\bar x}=d_2(\bar x_2v_2)$ and $u_2\in I_2(\bar y_2,v_2)$ (since
$u_2\in I_2(\bar y_2,s_2)$ and $v_2=q_2\in I_2(u_2,s_2)$).
Now $d_2(\bar x_2\bar y_2)=1$ implies $d_2(u_2v_2)=1$, as required.
  \end{proof}

For $i=0,\ldots,p$, define $P_i$ to be the subpath $x^i\ldots x^p$ of
$P$. Let $P_j$ be the maximal subpath with all edges in $E_2$ (i.e.,
$j$ is minimum subject to $x_1^j=\ldots=x_1^p$). Since $r(j)<r(j+1)<
\ldots <r(p)$, we can repeatedly apply Claim 1 to the edges of $P_j$,
starting with $x^{p-1}x^p$, and conclude that $N(x^i)$ is a singleton
$\{u^i\}$ with $u_1^i=y_1$ for each $i=j,\ldots,p$. Also $u^i=u^{i+1}$
if $\Delta_i<\Delta_{i+1}$, and $u^iu^{i+1}\in E_2$ if $\Delta_i=
\Delta_{i+1}$, where $\Delta_q$ stands for $\Delta^{x^q}$. Consider
two possible cases.

\medskip
{\em Case 1}\/\/: $j\ge 1$. By the maximality of $P_j$,
$x^{j-1}x^j\in E_1$. Let $b:=x_1^{j-1}$. For $i=j,\ldots, p$, define
$z^i$ and $v^i$ to be the points with $z^i_1=v^i_1=b$, $z^i_2=x^i_2$
and $v^i_2=u^i_2$, i.e., $z^i$ and $v^i$ are obtained by shifting the
points $x^i$ and $u^i$, respectively, along the edge $y_1b$ of $H_1$.
In particular, $z^j=x^{j-1}$. Denote $\Delta^{z^i}$  by $\Delta'_i$.

\bigskip
\noindent {\bf Claim 3.} {\em $\Delta'_i=\Delta_i$ and $v^i\in N(z^i)$
for each $i=j,\ldots,p$. }

\medskip
   \begin{proof}
Since $r(j-1)<r(j)$ and $x_2^{j-1}=x_2^j$, $r_1(b)<r_1(x^j)$. Therefore,
$u^i\in T$ implies $v^i\in T$, and we have $\Delta'_i\le d_2(z^iv^i)=
d_2(x^iu^i)=\Delta_i$. Suppose $\Delta'_i<\Delta_i$. Then $N(z^i)
\subseteq N(x^i)$, whence $N(z^i)=\{u^i\}$. But $d(z^iu^i)=d_1(by_1)
+d_2(x_2^iu_2^i)=1+d(z^iv^i)$; a contradiction. Thus, $\Delta'_i=
\Delta_i$ and $v^i\in N(z^i)$.
   \end{proof}

Consider the $x^{j-1}$--$y$ paths $P_{j-1}$ and $R=z^j\ldots z^px^p$
in $K$. From Claim 3 it follows that $\delta_1(z^iz^{i+1})=
\delta_1(x^ix^{i+1})$ for $i=j,\ldots,p-1$, and that
$\delta_1(x^{j-1}x^j)=\delta_1(z^px^p)$. Therefore,
$\delta_1(P_{j-1})=\delta_1(R)$. Replace in $P$ the part $P_{j-1}$
by $R$, forming the path $P'=x^0\ldots x^{j-1}z^j\ldots z^px^p
\ldots x^k$. Clearly $y=x^p$ is the first peak of $P'$. Cut off $y$
in $P'$ by replacing $y$ by a median $y''$ for $z^p,z,s$; let $P''$
be the resulting path. Since $z^py\in E_1$ and $yz\in E_2$, one
has $\rho (z^p,y,z,y'')\ge 0$, by Claim 1. Therefore, $\delta_1(P'')
\le\delta_1(P')=\delta_1(P)$, and~(\ref{eq:path_P}) follows by
induction because $z^p$ and $z$ are the first and second peaks of
$P''$ and $4^{r(y)}>4^{r(z^p)}+4^{r(z)}$.

\medskip
{\em Case 2}\/\/: $j=0$. Then $x^0=u^0$. By Claim 2 applied to the
edge $zy$, $N(z)$ is a singleton $\{\hat u\}$ with $\hat u_1=y_1$.
As before, let $y'\in M(x,z,s)\cap T$; then $y'_1=y_1$ and $N(y')$
is a singleton $\{v\}$ (by Claim 2 applied to the edge $y'x$).
Assuming $\Delta^{y'}<\Delta^x$ (equivalently: $\rho<0$), we have
$N(y')\subseteq N(x)\cap N(z)$. Hence, $v=\hat u=u^{p-1}$.

Form the $u^0$--$v$ path $R'$ by deleting repeated consecutive
elements in $u^0\ldots u^{p-1}$, and let $\bar R$ be the concatenation
of $R'$, a shortest $v$--$y'$ path $R''$, and the edge $y'x$.
Clearly the $\delta_1$-length of each edge of $R'$ is zero, while
the $\delta_1$-length of each edge of $R''$ is one. Also
$\delta_1(y'x)=1$.

Comparing $\bar R$ with the path $\bar P=x^0\ldots x^{p-1}$ and using
Claim 2, one can deduce that $\bar R=p-1$ (i.e., $\bar R$ is a
shortest path in $K$) and that $\delta_1(\bar R)=\delta_1(\bar P)$.
Now let $D$ be the concatenation of $R'$, $R''$ and the edge $y'z$.
Since $\delta_1(y'x)=\delta_1(y'z)$ and $\delta_1(xy)=\delta_1(yz)=0$,
we have  $\delta_1(D)=\delta_1(\bar R)
=\delta_1(P_0)$. Also $|D|=|R|=p-1$ implies that $D$
has no peaks. Then, replacing in $P$ the part $x^0\ldots x^{p+1}$ by
$D$, we obtain
the path $P'$ with $\delta_1(P')=\delta(P)$ and $\omega(P')<
\omega(P)$ and can apply induction.

Thus, (\ref{eq:path_P}) is proven. In order to conclude that $m_1$ is an
extension of $d_1$, it suffices to consider a path $L$ as
in~(\ref{eq:path_P}) and show the following:
   \begin{myitem}
 (i) $\;\;$ if $z^\alpha\in T$, then $\delta_1(L)\ge
d_1(z_1^0z_1^\alpha)$;
   \begin{itemize}
 \item[(ii)] $\delta_1(L)+\delta_1(z^\alpha s_1)\ge r_1(z_1^0)$.
    \end{itemize}
   \label{eq:z_alpha}
   \end{myitem}
(In fact, (i) embraces the case of a path in $G_1$, with both ends
in $T_1$, whose first and last edges belong to $\Ascr_1$, while (ii)
does the case when one of these edges is in $\Ascr_1$ and the other
in $\Bscr_1$.) Case (i) is trivial because $z^\alpha\in T$ means
that $L$ is a path in $H$, and therefore, the $\delta_1$-length of
each of its edges in $E_1$ is equal to one. So let us prove (ii).
One may assume that $r(z^0)<\ldots<r(z^\alpha)$ (taking into account
that $d_1(z_1^0s_1)\le d_1(z_1^0z_1^\beta)+d_1(z_1^\beta z_1^\alpha)$
and $\delta_1(L)=\delta_1(L')+\delta_1(L'')$, where $L'=z^0\ldots
z^\beta$ and $L''=z^\beta\ldots z^\alpha$, and assuming w.l.o.g.
that $L'$ is $\delta_1$-shortest).

For $i=0,\ldots,\alpha$, let $\ell_i$ denote the $\delta_1$-length
of the path $z^0\ldots z^i$, and let $\rho_i$ and $\Delta_i$ stand
for $r_1(z_1^i)$ and $\Delta^{z^i}$, respectively. By the
definition of $\delta_1$ on $\Bscr_1$, $\delta_1(z^is_1)$ is equal
to $\rho_i-\Delta_i$. We show that
   \begin{equation} \label{eq:ell_rho}
  \ell_i+\rho-\Delta_i\ge \rho_0,
    \end{equation}
using induction on $i$. This gives the desired
inequality~(\ref{eq:z_alpha})(ii) when $i=\alpha$. Since $\ell_0=
\Delta_0=0$, (\ref{eq:ell_rho}) holds for $i=0$. Assume it holds
for $i-1$ ($0<i<\alpha$), and let $a:=\ell_i-\ell_{i-1}$, $b:=
\rho_i-\rho_{i-1}$ and $c:=\Delta_i-\Delta_{i-1}$.
Then~(\ref{eq:ell_rho}) for $i$ follows from $a+b-c\ge 0$. To see
the latter, consider four possible cases for $e=z^{i-1}z^i$,
taking into account that $\Delta_i\ge \Delta_{i-1}$ since
$r(z^i)>r(z^{i-1})$.

(a) Let $e\in E_1$ and $\Delta_i=\Delta_{i-1}$. Then $a+b-c=
1+1+0=2$.

(b) Let $e\in E_1$ and $\Delta_i>\Delta_{i-1}$. Then $a+b-c=
0+1-1=0$.

(c) Let $e\in E_2$ and $\Delta_i=\Delta_{i-1}$. Then $a+b-c=
0+0-0=0$.

(d) Let $e\in E_2$ and $\Delta_i>\Delta_{i-1}$. Then $a+b-c=
1+0-1=0$.

Thus, $m_1$ is an extension of $d_1$. It remains to show that
$m_i(e)=\delta_i(e)$ for $i=1,2$ and $e\in\Escr_i$. This is
obvious when $e\in U\cup U_i$ or when $\delta_i(e)=0$. If
$e=xs_i\in\Bscr_i$, then $m_i(e)=\delta_i(e)$ follows from the
fact that for $t\in N(x)$, the path in $G_i$ obtained by
concatenationg the edge $t_it$, a shortest $t$--$x$ path in $K$,
and the edge $xs_i$ is $\delta_i$-shortest (this fact was shown
at the beginning of this section). Finally, each edge $e\in E(K)$
belongs to a shortest $t$--$t'$ path $P$ in $K$ with $t,t'\in T$.
Since $\delta_1(e')+\delta_2(e')=1$ for all edges $e'$ of $K$,
we have $\delta_1(P)+\delta_2(P)=|P|=d(tt')=d_1(t_1t'_1)+
d_2(t_2t'_2)$, whence $\delta_i(P)=d_i(t_it'_i)$, implying
$m_i(e)=\delta_i(e)$.

This completes the proof of Lemma~\ref{lm:m_i} and completes the
proof of Theorems~\ref{tm:gen_retrac} and~\ref{tm:main}.

%
%
\section{Intractable Cases} \label{sec:nphard}

In this section we prove Theorem~\ref{tm:nphard}, considering a
metric $\mu$ on a set $T$ such that either $\mu$ is non-modular or
$\mu$ is modular but its underlying graph $H=(T,U)$ is
non-orientable. W.l.o.g., one may assume $\mu$ is integer-valued.
Our method borrows the idea from~\cite{kar-98} for the path metrics
$\mu=d^H$ as in Theorem~\ref{tm:hardgr}, which in turn generalizes
the construction from~\cite{DJ-94} for $H=K_3$.

Given a set $V\supset T$, a function $E_V\to\Zset_+$, nodes $s,t\in T$,
and points $x,y\in V-T$, let $\tau(s,x|t,y)$ denote the minimum
$c\cdot m$ among all $m\in \exten^0(\mu,V)$ such that $m(xs)=m(yt)=0$.

The core of the proof in~\cite{DJ-94} that the 3-terminal cut problem
is NP-hard is the construction of a ``gadget'' $(V,c)$ with
specified $s,t,x,y$ satisfying the following property:
   \begin{myitem}
  (i)$\;\;$ $\tau(s,x|t,y)=\tau(s,y|t,x)=\hat\tau$,
     \begin{itemize}
  \item[(ii)]$\;\;$ $\tau(s,x|s,y)=\tau(t,x|t,y)=\hat\tau+\delta$
for some $\delta>0$,
  \item[(iii)]$\;\;$ $\tau(s',x|t',y)\ge\hat\tau+\delta$ for all
other pairs $\{s',t'\}$ in $T$,
     \end{itemize}
   \label{eq:sxty}
   \end{myitem}
where $\hat\tau$ stands for $\tau(V,c,\mu)$ (with $\mu=d^{K_3}$).
Then the NP-hardness of the problem is easily shown by a reduction
from MAX CUT.

Our aim is to construct corresponding ``gadgets'' satisfying
(\ref{eq:sxty}) for $\mu$ as in Theorem~\ref{tm:nphard}; then the
theorem will follow by a similar reduction.

First we consider the case when $\mu$ is modular but $H$ is
non-orientable, which is technically simpler. In fact, the
construction and arguments in this case are similar to those for the
corresponding unweighted case ($\mu=d^H$) given
in~\cite[Sec. 6]{kar-98}. More precisely, since $H$ is non-orientable,
there exists a projective sequence $(e_0,e_1,\ldots,e_{k-1},e_k=e_0)$
of edges of $H$ yielding the ``twist'' (or forming the
{\em orientation-reversing dual cycle}). That is,
   \begin{myitem}
  for $i=0,\ldots,k-1$, $e_i=s_it_i$ and $e_{i+1}=s_{i+1}t_{i+1}$ are
opposite edges in the 4-circuit $C_i=s_it_it_{i+1}s_{i+1}s_i$, and
$t_k=s_0$ (and $s_k=t_0$).
    \label{eq:dual_cycle}
   \end{myitem}
(One can choose such a sequence with all edges (though not necessarily
the nodes) distinct, but this is not important for us.) Since $\mu$
is modular, we have by~(\ref{eq:mod2}) that
   \begin{myitem}
 for $i=0,\ldots,k-1$, $\mu(e_i)$ is a constant $h$, and
$\mu(s_is_{i+1})=\mu(t_it_{i+1})=:f_i$.
   \label{eq:h_and_f}
   \end{myitem}

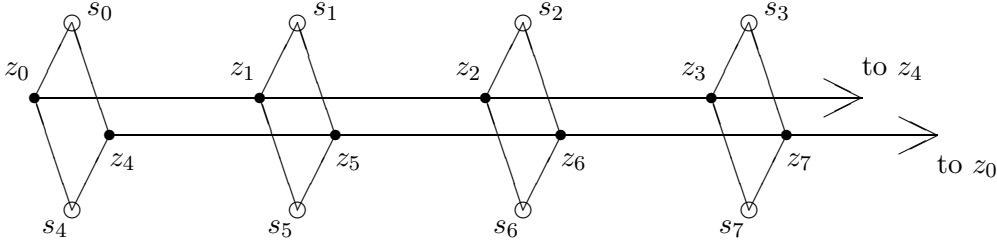
\begin{figure}[tb]
 \unitlength=1.0mm
\begin{center}
  \begin{picture}(135,35)
\multiput(10,5)(30,0){4}%
{\put(0,0){\circle{2}}
\put(0,25){\circle{2}}
\put(-5,15){\circle*{1.5}}
\put(5,10){\circle*{1.5}}
\put(0,0){\line(-1,3){5}}
\put(5,10){\line(-1,3){5}}
\put(0,0){\line(1,2){5}}
\put(-5,15){\line(1,2){5}}
}
\put(5,20){\line(1,0){110}}
\put(15,15){\line(1,0){110}}
\put(115,20){\line(-2,1){5}}
\put(125,15){\line(-2,1){5}}
\put(115,20){\line(-2,-1){5}}
\put(125,15){\line(-2,-1){5}}
\put(12,31){$s_0$}
\put(42,31){$s_1$}
\put(72,31){$s_2$}
\put(102,31){$s_3$}
\put(6,2){$s_4$}
\put(36,2){$s_5$}
\put(66,2){$s_6$}
\put(96,2){$s_7$}
\put(1,23){$z_0$}
\put(31,23){$z_1$}
\put(61,23){$z_2$}
\put(91,23){$z_3$}
\put(15,11){$z_4$}
\put(45,11){$z_5$}
\put(75,11){$z_6$}
\put(105,11){$z_7$}
\put(115,23){to $z_4$}
\put(125,10){to $z_0$}
  \end{picture}
 \end{center}
\caption{\hspace{2cm} gadget for a non-orientable $H$ \hspace{2cm}}
\label{fig:non_orient}
  \end{figure}

We denote $t_i$ by $s_{i+k}$ and take indices modulo $2k$. The
desired gadget is represented by the graph $G=(V,E)$ with the weights
$c(e)$ of edges $e\in E$, where $V=T\cup\{z_0,\ldots,z_{2k-1}\}$ and for
$i=0,\ldots,2k-1$,
  \begin{itemize}
 \item[(i)] $z_i$ is adjacent to both $s_i$ and $s_{i+k}$, and $c(z_is_i)=
c(z_is_{i+k})=N$ for a positive integer $N$ (specified below);
 \item[(ii)] $z_i$ and $z_{i+1}$ are adjacent, and $c(z_iz_{i+1})=1$.
   \end{itemize}
Figure~\ref{fig:non_orient} illustrates $G$ for $k=4$. We put $s=s_0$,
$t=t_0$, $x=z_0$ and $y=z_k$, and formally extend $c$ by zero to $E_V-E$.
We assert that~(\ref{eq:sxty}) holds.
%
%

Indeed, each $m\in\exten^0(\mu,V)$ is associated with the mapping
$\gamma:\{z_0,\ldots,z_{2k-1}\}\to T$, where $\gamma(z_i)=s_j$ if
$m(z_is_j)=0$; we say that $z_i$ is {\em attached}\/ by $\gamma$
to $s_j$ and denote $m$ by $m^\gamma$. If $\gamma(z_i)=v$, then,
letting $\eps:=\mu(s_iv)+\mu(vs_{i+k})-\mu(s_is_{i+k})$, the
contribution to the {\em volume}\/ $c\cdot m^\gamma$ due to the edges
$e=z_is_i$ and $e'=z_is_{i+k}$ is equal to
  $$
 c(e)m^\gamma(e)+c(e')m^\gamma(e')=N(m^\gamma(e)+m^\gamma(e'))
   =Nh+N\eps ;
  $$
cf. (\ref{eq:h_and_f}). We have $\eps=0$ if $v\in\{s_i,s_{i+k}\}$,
and $\eps\ge 1$ otherwise. Hence, every mapping $\gamma$ pretending
to be optimal or nearly optimal must attach each $z_i$ to either
$s_i$ or $s_{i+k}$ whenever $N$ is chosen sufficiently
large (e.g., $N=1+2k\max\{\mu(st): s,t\in T\}$).

Next, if $z_i$ is attached to $s_i$ (resp. $s_{i+k}$) and $z_{i+1}$
to $s_{i+1}$ (resp. $s_{i+1+k}$), then the edge $u=z_iz_{i+1}$
contributes $c(u)m^\gamma(u)=f_i$ (cf.~(\ref{eq:h_and_f}),
letting $f_{j}=f_{j+k}$. On the other
hand, if $z_i$ is attached to $s_i$ (resp. $s_{i+k}$) while $z_{i+1}$
to $s_{i+1+k}$ (resp. $s_{i+1}$), then the contribution becomes $h+f_i$
($=\mu(s_it_{i+1})$).

So we can conclude that $\hat\tau=
2khN+2(f_1+\ldots+f_k)$, and there are precisely two optimal
0-extensions, namely, $m^{\gamma_1}$ and $m^{\gamma_2}$,
where $\gamma_1(z_i)=s_i$ and $\gamma_2(z_i)=s_{i+k}$ for $i=0,\ldots,
2k-1$. This gives (i) in~(\ref{eq:sxty}).
Furthermore, one can see that if $m^\gamma$ is the least-volume
0-extension induced by $\gamma$ that brings both $x,y$ either to
$s$ or to $t$, then $m^\gamma(z_jz_{j+1})=h+f_j$ for precisely
two numbers $j\in\{0,\ldots,2k-1\}$ such that
$f_j=\min\{f_1,\ldots,f_k\}$. So $c\cdot m^\gamma=\hat\tau+2h$,
yielding~(\ref{eq:sxty})(ii). Finally, (iii) is ensured by the choice
of $N$.

Thus, (\ref{eq:0ext}) with $\mu$ modular and $H$ non-orientable
is NP-hard. Moreover, it is strongly NP-hard because the number
$N$ is a constant depending only on $\mu$.

Next we consider the case when $\mu$ is not modular. Let
$\Delta(x,y,z)$ denote the value ({\em perimeter}) $\mu(xy)+\mu(yz)
+\mu(zx)$ for $x,y,z\in T$. We fix a medianless triplet
$\{s_0,s_1,s_2\}$ such that $\Delta(s_0,s_1,s_2):=\bar\Delta$
is minimum. By
technical reasons, we put $s_{i+3}=s_i$, $i=0,1,2$, and take indices
modulo 6. The gadget $(G=(V,E),c)$ that we construct has a somewhat
more complicated structure compared with that for the corresponding
unweighted case in~\cite[Sec. 6]{kar-98}. Here
   $$
 V=T\cup Z,\quad Z=\{z_0,\ldots,z_5\}\quad\mbox{and}\quad
 E=E_1\cup E_2\cup E_3.
   $$
For $i=1,2,3$, the edges $e\in E_i$ are endowed with weights
$c_i(e)$, and $c(e)$ is defined to be $N_ic_i(e)$. The
factors $N_1,N_2,N_3$ are chosen so that $N_1=1$, $N_2$ is
sufficiently large, and $N_3$ is sufficiently large with respect
to $N_2$. Informally speaking, the ``heavy'' edges of $E_3$ provide
that (at optimality or almost optimality) each point $z_j$ gets into
the interval $I_j:=\{v\in T:\mu(s_{j-1}v)+\mu(vs_{j+1})=
\mu(s_{j-1}s_{j+1})\}$, then the ``medium'' edges of $E_2$ make $z_j$
choose only between the endpoints $s_{j-1},s_{j+1}$ of $I_j$, and
finally the ``light'' edges of $E_1$ provide the desired property
(\ref{eq:sxty}).

As before, $m^\gamma$ denotes the 0-extension of $\mu$ to $V$ induced
by $\gamma:Z\to T$. Define $d_i:=d_{i+3}=\mu(s_{i-1}s_{i+1})$.
We say that a path $P=(v_1,\ldots,v_k)$ on $T$ is shortest if it
is $\mu$-shortest.

The set $E_3$ consists of the edges $e_j=z_js_{j-1}$
and $e'_j=z_js_{j+1}$ with $c_3(e_j)=c_3(e'_j)=1$ for $j=0,\ldots,5$.
Then the contribution to $c\cdot m^\gamma$ due to $e_j$ and $e'_j$
is $N_3d_j$ if $\gamma(z_j)\in I_j$, and at least $N_3d_j+N_3$
otherwise, yielding that $z_j$ should be mapped into $I_j$, by the
choice of $N_3$. The minimality of $\bar\Delta$ provides the
following useful property.
%
%
    \begin{statement} \label{st:about_v}
For any $v\in I_j$, at least one of the paths $P=(s_j,s_{j-1},v)$
and $P'=(s_j,s_{j+1},v)$ is shortest.
    \end{statement}
   \begin{proof}
Let for definiteness $j=1$. Suppose $P'$ is not shortest. Then
$\mu(s_1v)<|P'|=\mu(s_1s_2)+\mu(s_2v)$ and
$\mu(s_0v)=\mu(s_0s_2)-\mu(s_2v)$
imply $\Delta(s_1,v,s_0)<\bar\Delta$. So $s_1,v,s_0$
have a median $w$. If $w=s_0$, $P$ is shortest. Otherwise we have
$\Delta(s_1,w,s_2)<\Delta$ (since $\mu(s_1w)<\mu(s_1s_0)$ and the
path $(s_2,v,w,s_0)$ is, obviously, shortest). Then $s_1,w,s_2$
have a median $q$. It is easy to see that $q$ is a median for
$s_0,s_1,s_2$; a contradiction.
   \end{proof}

We now explain the construction of $E_2$ and $c_2$. Each $z=z_j$
($j=0,\ldots,5$) is connected to each $s_i$ ($i=0,1,2$) by edge
$u_i=zs_i$ whose weight is defined by
   \begin{equation} \label{eq:a_i}
  c_2(u_i)=(d_{i-1}+d_{i+1}-d_i)/(d_{i-1}d_{i+1})=:a_i
   \end{equation}
($a_i$ is positive and does not depend on $j$).
Suppose $z$ is mapped by $\gamma$
to some $s_i$, say $\gamma(z)=s_1$. Then, up to a factor of $N_2$,
the contribution to $c\cdot m^\gamma$ from the edges $u_0,u_1,u_2$
(concerning $z$) is
     \begin{eqnarray}
   && d_2a_0+d_0a_2=d_2(d_1+d_2-d_0)/(d_1d_2)+
            d_0(d_1+d_0-d_2)/(d_0d_1)  \label{eq:dada}\\
     && \qquad\qquad\qquad (d_1+d_2-d_0)/d_1+(d_1+d_0-d_2)/d_1=2.
                                     \nonumber
  \end{eqnarray}
On the other hand, the contribution grows when $z_j$ falls into
the interior of any interval $I_i$.
%
%
   \begin{statement} \label{st:interior}
 Let $v\in I_i-\{s_{i-1},s_{i+1}\}$. Then $\sigma:=
\sum(a_i\mu(s_iv): i=0,1,2)>2$.
   \end{statement}
  \begin{proof}
Let for definiteness $i=0$, $\mu(s_1v)=\eps$ and $\mu(s_0v)=d_2+\eps$
(cf. Statement~\ref{st:about_v}). Then
   $$
 \sigma=(d_2+\eps)a_0+\eps a_1+(d_0-\eps)a_2=d_2a_0+d_0a_2+
\eps(a_0+a_1-a_2)=2+\eps(a_0+a_1-a_2),
   $$
in view of~(\ref{eq:dada}). We observe that $a_0+a_1-a_2>0$. Indeed,
   $$
  \begin{array}{l}
 \quad d_0d_1d_2(a_0+a_1-a_2)=(d_0d_1+d_0d_2-d_0^2)+
    (d_1d_0+d_1d_2-d_1^2)-(d_2d_0+d_2d_1-d_2^2) \\
  \hfill =2d_0d_1-d_0^2-d_1^2+d_2^2=d_2^2-(d_0-d_1)^2>0
   \end{array}
  $$
since $d_2>d_0-d_1$. So $\sigma>2$.
   \end{proof}

Thus, by an appropriate choice of constants $N_2$ and $N_3$,
each point $z_j$ must be mapped to either $s_{j-1}$ or $s_{j+1}$.
Such a mapping $\gamma$ is called {\em feasible}.
We now construct the crucial
set $E_1$ and function $c_1$. The set $E_1$ consists of six edges
$g_j=z_jz_{j+1}$, $j=0,\ldots,5$, forming the 6-circuit $C$ (this
is similar to the construction in~\cite{kar-98} motivated
by~\cite{DJ-94}). The essense is how to assign $c_1$. For $i=0,1,2$,
let $h_i:=h_{i+3}:=(d_{i-1}+d_{i+1}-d_i)/2$. These numbers would
be just the distances from $s_0,s_1,s_2$ to their median if it
existed, i.e.,
   \begin{equation} \label{eq:d_i}
  d_i=h_{i-1}+h_{i+1}.
   \end{equation}
We define
   \begin{equation}  \label{eq:c_one}
 c_1(z_jz_{j+1})=c_1(z_{j+3}z_{j+4})=h_{j-1} \qquad\mbox{for}
                                       \quad j=0,1,2.
    \end{equation}

For $\gamma:Z\to T$, let $\zeta^\gamma$ denotes
$\sum(c_1(g_j)m^\gamma(g_j): j=0,\ldots,5)$, i.e., $\zeta^\gamma$ is
the contribution to $c\cdot m^\gamma$ from the edges of $C$. The
analysis below will depend on the numbers
    \begin{equation}  \label{eq:rho_alpha}
  \rho=2(h_0h_1+h_1h_2+h_2h_0)\quad\mbox{and}\quad
           \alpha=2\/\/\/ \min\{h_0^2,h_1^2,h_2^2\}.
    \end{equation}
W.l.o.g., assume $h_0\le h_1,h_2$, i.e., $2h_0^2=\alpha$.
Our aim is to show that
(\ref{eq:sxty}) holds if we take as $s,t,x,y$ the elements
$s_0,s_2,z_1,z_4$, respectively.

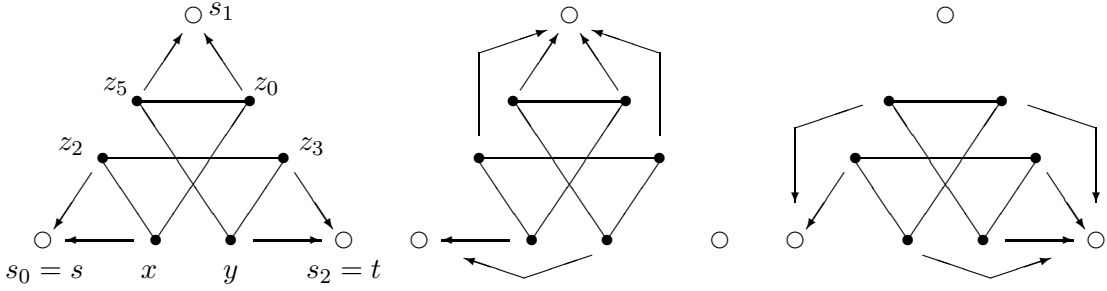
\begin{figure}[tb]
 \unitlength=1.0mm
\begin{center}
  \begin{picture}(150,35)
\multiput(0,0)(50,0){3}%
{\put(5,5){\circle{2}}
\put(45,5){\circle{2}}
\put(25,35){\circle{2}}
\put(20,5){\circle*{1.5}}
\put(30,5){\circle*{1.5}}
\put(13,16){\circle*{1.5}}
\put(37,16){\circle*{1.5}}
\put(17.5,23.5){\circle*{1.5}}
\put(32.5,23.5){\circle*{1.5}}
\put(13,16){\line(1,0){24}}
\put(17.5,23.5){\line(1,0){15}}
\put(20,5){\line(-2,3){7}}
\put(20,5){\line(2,3){12.5}}
\put(30,5){\line(-2,3){12.5}}
\put(30,5){\line(2,3){7}}
}
\put(17,5){\vector(-1,0){9}}
\put(67,5){\vector(-1,0){9}}
\put(33,5){\vector(1,0){9}}
\put(133,5){\vector(1,0){9}}
\put(11.5,14){\vector(-2,-3){5}}
\put(111.5,14){\vector(-2,-3){5}}
\put(38.5,14){\vector(2,-3){5}}
\put(138.5,14){\vector(2,-3){5}}
\put(18.5,25){\vector(2,3){5}}
\put(68.5,25){\vector(2,3){5}}
\put(31.5,25){\vector(-2,3){5}}
\put(81.5,25){\vector(-2,3){5}}
\put(69,0){\line(3,1){9}}
\put(69,0){\vector(-3,1){8}}
\put(131,0){\line(-3,1){9}}
\put(131,0){\vector(3,1){8}}
\put(63,30){\line(0,-1){11}}
\put(63,30){\vector(3,1){9}}
\put(87,30){\line(0,-1){11}}
\put(87,30){\vector(-3,1){9}}
\put(105,20){\line(3,1){9}}
\put(105,20){\vector(0,-1){10}}
\put(145,20){\line(-3,1){9}}
\put(145,20){\vector(0,-1){10}}
\put(0,0){$s_0=s$}
\put(18,0){$x$}
\put(29,0){$y$}
\put(40,0){$s_2=t$}
\put(7,17){$z_2$}
\put(39,17){$z_3$}
\put(13,25){$z_5$}
\put(33,25){$z_0$}
\put(27,35){$s_1$}
  \end{picture}
 \end{center}
\caption{\hspace{0.5cm} (a) $\gamma_1$ \hspace{3.5cm} (b) $\gamma_3$
    \hspace{3.5cm} (c) $\gamma_4$ \hspace{1cm}}
\label{fig:gadget}
  \end{figure}

To show this, consider the mapping $\gamma_1$ as drawn in
Fig.~\ref{fig:gadget}a, i.e., $\gamma_1(z_j)$ is $s_{j+1}$ for
$j=0,2,4$ and $s_{j-1}$ for $j=1,3,5$. This $\gamma_1$ attaches
$x$ to $s$ and $y$ to $t$. In view
of~(\ref{eq:d_i})--(\ref{eq:rho_alpha}), we have
   $$
    \begin{array}{l}
 \quad \zeta^{\gamma_1}=c_1(g_0)\mu(\gamma_1(z_0)\gamma_1(z_1))+\ldots+
   c_1(g_5)\mu(\gamma_1(z_5)\gamma_1(z_0)) \hfill\\
  \hfill =h_2d_2+h_0\cdot 0+h_1d_1+h_2\cdot 0+h_0d_0+h_1\cdot 0
    =h_2(h_0+h_1)+h_1(h_0+h_2)+h_0(h_1+h_2)=\rho. \quad
    \end{array}
   $$
Similarly, $\zeta^{\gamma_2}=\rho$ for the symmetric mapping
$\gamma_2$ which is defined by $\gamma_2(z_j)=\gamma_1(z_{j+3})$,
attaching $x$ to $t$ and $y$ to $s$. We shall see later that
$\gamma_1$ and $\gamma_2$ are just optimal mappings for our
gadget.
%
%

The mappings pretending to provide (ii) in~(\ref{eq:sxty}) are
$\gamma_3$ and $\gamma_4$ illustrated in Fig.~\ref{fig:gadget}b,c;
here both $x,y$ are mapped by $\gamma_3$ to $s$, and by $\gamma_4$
to $t$. We have
   $$
    \begin{array}{l}
 \quad \zeta^{\gamma_3}=h_2d_2+h_0d_2+h_1\cdot 0+h_2d_2+h_0d_2+
          h_1\cdot 0 =(2h_2+2h_0)(h_0+h_1) \hfill \\
  \hfill =2h_2h_0+2h_2h_1+2h_0^2+2h_0h_1=\rho+\alpha \quad
    \end{array}
   $$
and
   $$
    \begin{array}{l}
 \quad \zeta^{\gamma_4}=h_2\cdot 0+h_0d_1+h_1d_1+h_2\cdot 0+
   h_0d_1+h_1d_1=(2h_0+2h_1)(h_0+h_2) \hfill  \\
  \hfill =2h_0^2+2h_0h_2+2h_1h_0+2h_1h_2=\rho+\alpha. \quad
    \end{array}
   $$

Now (\ref{eq:sxty}) is implied by the following.
%
%
   \begin{statement}  \label{st:zeta_gamma}
Let $\gamma$ be a feasible mapping different from $\gamma_1$ and
$\gamma_2$. Then $\zeta^\gamma\ge\rho+\alpha$.
    \end{statement}
   \begin{proof}
By~(\ref{eq:d_i}), $\zeta^\gamma$ is representable as a nonnegative
integer combination of products $h_ih_j$ for $0\le i,j\le 5$
(including $i=j$). The contribution $\zeta_j$ to $c\cdot m^\gamma$
from a single edge $g_j=z_jz_{j+1}$ is as follows:
     \begin{myitem}
  (i)$\;\;$ if $\gamma(z_j)=\gamma(z_{j+1})=s_{j-1}$, then
$\zeta_j=0$;
     \begin{itemize}
\item[(ii)]
  if $\gamma(z_j)=s_{j+1}$ and $\gamma(z_{j+1})=s_j$, then $\zeta_j=
h_{j-1}d_{j-1}=h_{j-1}h_j+h_{j-1}h_{j+1}$;
\item[(iii)]
  if $\gamma(z_j)=s_{j+1}$ and $\gamma(z_{j+1})=s_{j-1}$, then $\zeta_j=
h_{j-1}d_{j}=h_{j-1}h_{j+1}+h_{j-1}^2$;
\item[(iv)]
  if $\gamma(z_j)=s_{j-1}$ and $\gamma(z_{j+1})=s_{j}$, then $\zeta_j=
h_{j-1}d_{j+1}=h_{j-1}h_{j}+h_{j-1}^2$;
   \end{itemize}
     \label{eq:zeta_j}
     \end{myitem}

We call $g_j$ {\em slanting}\/ if it is as in case (iii) or (iv)
of~(\ref{eq:zeta_j}). If no edge of $C$ is slanting, then $\gamma$
is either $\gamma_1$ or $\gamma_2$. Otherwise $C$ contains at
least two slanting edges. In this case we observe
from~(\ref{eq:zeta_j}) that the representation of $\zeta^\gamma$
includes $h_i^2+h_j^2$ (or $2h_i^2$) for some $i,j$, which is at
least $\alpha$. Now the result follows from the fact that the
representation includes $2h_ih_j$ for each $0\le i<j\le 2$.

To see the latter, w.l.o.g., assume $i=0$, $j=2$, and consider the
edges $g_0$ and $g_1$. By~(\ref{eq:d_i}), $g_0$ contributes $h_0h_2$
in cases (ii),(iv), i.e., when $\gamma(z_1)=s_0$. And if
$\gamma(z_1)=s_2$, then $g_1$ contributes $h_0h_2$. Similarly, the
pair $g_3,g_4$ contributes $h_0h_2$.
   \end{proof}

This completes the proof of Theorem~\ref{tm:nphard}.

\Xcomment{
\section{Concluding Remarks} \label{sec:concl}
}

\small
\bibliographystyle{plain}
\bibliography{uncross}

\end{document}